\documentclass[11pt]{article}

\usepackage{authblk}
\usepackage{amsfonts}
\usepackage{amsthm}
\usepackage{amsmath}
\usepackage{graphicx}
\usepackage{comment}
\usepackage[usenames,dvipsnames]{xcolor}
\usepackage{enumerate}
\usepackage{extarrows}
\usepackage{hyperref}
\usepackage[round, authoryear]{natbib}
\usepackage{stackrel}
\usepackage{centernot}
\usepackage{caption}
\usepackage{subcaption}
\usepackage{fancyhdr}

\let\quoteOLD\quote
\def\quote{\quoteOLD\small}

\definecolor{labelkey}{cmyk}{0,0.8,1,0.5}
\definecolor{refkey}{cmyk}{0,0.8,1,0.5}

\newtheorem{theorem}{Theorem}[section]
\newtheorem{example0}{\sc Example}[subsection]

\newtheorem{corollary}{Corollary}
\newtheorem{proposition}{Proposition}

\newtheorem{remark}{Remark}

\numberwithin{equation}{section}
\numberwithin{theorem}{section}
\numberwithin{corollary}{section}
\numberwithin{proposition}{section}
\numberwithin{lemma}{section}
\numberwithin{definition}{section}
\numberwithin{remark}{section}

\makeatletter
\def\th@newremark{\th@remark\thm@headfont{\bfseries}}
\makeatletter

\def\boxit#1{\vbox{\hrule\hbox{\vrule\kern6pt
          \vbox{\kern6pt#1\kern6pt}\kern6pt\vrule}\hrule}}

\newcommand{\pibar}{\overline{\Pi}}
\newcommand{\pibarinv}{\overline{\Pi}^{\leftarrow}}
\newcommand{\pibarpinv}{\overline{\Pi}^{+,\leftarrow}}

\newcommand{\lambar}{\overline{\Lambda}}

\newcommand{\veps}{\varepsilon}

\newcommand{\topr}{\stackrel{\mathrm{P}}{\longrightarrow}}
\newcommand{\todr}{\stackrel{\mathrm{D}}{\longrightarrow}}
\newcommand{\eqdr}{\stackrel{\mathrm{D}}{=}}

\newcommand{\toweak}{\stackrel{w}{\longrightarrow}}

\newcommand{\R}{\Bbb{R}}

\newcommand{\N}{\Bbb{N}}

\newcommand{\rmd}{{\rm d}}
\newcommand{\rmi}{{\rm i}}
\newcommand{\halmos}{\quad\hfill\mbox{$\Box$}}

\newcommand{\wt}{\widetilde}

\newcommand{\dto}{\downarrow}

\newcommand{\SSSS}{\Gamma}

\newcommand{\be}{\begin{equation}}
\newcommand{\ee}{\end{equation}}
\newcommand{\bea}{\begin{eqnarray}}
\newcommand{\eea}{\end{eqnarray}}
\newcommand{\bean}{\begin{eqnarray*}}
\newcommand{\eean}{\end{eqnarray*}}
\newcommand{\ben}{\begin{equation*}}
\newcommand{\een}{\end{equation*}}
\newcommand{\ba}{\begin{aligned}}
\newcommand{\ea}{\end{aligned}}

\def\nexto{\kern -0.54em}

\def\topr{\buildrel P \over \to }

\newcommand{\PP}{\textbf{\rm P}}
\newcommand{\EE}{\textbf{\rm E}}
\newcommand{\PD}{\textbf{\rm PD}}

\newcommand{\bfx}{{\bf x}}

\begin{document}

\title{\bf Convergence to Stable Limits for Ratios of Trimmed L\'evy Processes and their Jumps}
\author{Yuguang F. Ipsen,  P\'eter Kevei\thanks{Research supported by the J\'anos Bolyai Research Scholarship of the Hungarian 
Academy of Sciences, and by NKFIH grant FK124141. }
\  and Ross A. Maller\thanks{Research
partially supported by ARC Grant DP1092502. \newline
Email: Yuguang.Ipsen@anu.edu.au; kevei@math.u-szeged.hu; Ross.Maller@anu.edu.au}}

\affil{Research School of Finance, Actuarial Studies \& Statistics\\
Australian National University
 \\Canberra, ACT, 0200, Australia\\
and \\
MTA-SZTE Analysis and Stochastics Research Group\\
Bolyai Institute, Aradi v\'ertan\'uk tere 1\\
 6720 Szeged, Hungary
}

\maketitle

\begin{abstract}
We derive characteristic function  identities for conditional distributions of an $r$-trimmed L\'evy process given  its  $r$ largest jumps up to a designated time $t$. 
Assuming the underlying L\'evy process is in the domain of attraction of a stable process as $t\dto 0$, these identities  are applied to show joint convergence of the trimmed process  divided by  its large jumps to corresponding quantities constructed from a stable limiting process.
This generalises related results in the 1-dimensional subordinator case developed in
\cite{KeveiMason2014} and produces new discrete distributions on the infinite simplex in the limit.
\end{abstract}

\section{Introduction and L\'evy Process Setup}\label{intro}
Deleting the $r$ largest jumps up to a designated time $t$ 
from a L\'evy process gives the ``$r$-trimmed L\'evy process". 
We derive useful characteristic function  identities for conditional distributions of the process given  some of its  largest jumps.
 As corollaries, representations for the characteristic functions of the trimmed process  divided by  its large jumps are found. Assuming $X$ is in the domain of attraction of a stable process as $t\dto 0$, the representations  are applied to show  joint convergence of those ratios 
to corresponding quantities constructed from  the stable limiting process.

In the case of subordinators,  \cite{KeveiMason2014} considered one-dimensional convergence to stable subordinators and derived the limit distribution of the ratio of an $r$-trimmed subordinator to its $r^{th}$ largest jump occurring up till a specified time $t>0$, as $t\dto 0$ or $t\to \infty$. 
\cite{perman1993}, also considering subordinators, derived exact expressions for the joint density of the ratios of the first $r$  largest  jumps up till time $t=1$ of a subordinator, taken as ratios of the value of the subordinator itself at time 1. In Perman's case the canonical measure of the subordinator was assumed to have a density with respect to Lebesgue measure. His results, when applied to a Gamma subordinator, produce formulae for the  Poisson-Dirichlet process.

For our asymptotic results we allow a general L\'evy measure, making no continuity assumptions on it.
Our main result, Theorem \ref{stable}, 
is a multivariate version of part of Theorem 1.1 of \cite{KeveiMason2014}, and,  as a generalisation, we consider a trimmed 
L\'evy process in the domain of attraction of  a stable distribution with parameter $\alpha$ in $(0,2)$,
taken as a ratio of one of its large jumps at time $t$. We show the joint convergence of these ratios   to corresponding quantities constructed from  the stable limiting process, as time $t$ tends to 0.

When $0<\alpha <1$, the limit distribution in Theorem \ref{stable} is related to
the generalised Poisson-Dirichlet distribution $\PD_\alpha^{(r)}$ in \cite{trimfrenzy} derived from the trimmed stable subordinator, which includes as a special case the $\PD(\alpha, 0)$ distribution in \cite{PY1997}. 
When $\alpha>1$ the process is not a subordinator, and there is no direct connection with the Poisson-Dirichlet distribution.
In this case the process has  to be centered
appropriately to get the required convergence.
We note that (since the L\'evy measure has infinite mass) there are always infinitely many ``large" jumps of $X_t$, a.s., in any 
right neighbourhood of 0.

These considerations form the basis of further generalised versions of  Poisson-Dirichlet distributions explored in \cite{trimfrenzy}.
In the present paper we limit ourselves to proving Theorem \ref{stable} (in Section \ref{s5}) and the foundational results needed for its proof (in Section \ref{s3}).
A second  Theorem \ref{prin} proves a kind of ``large trimming" result, showing that the trimmed process is of small order of the largest jump trimmed, uniformly in $t$, as the order tends to infinity.
Section \ref{ps5} contains the proofs of the results in Section \ref{s5}.
For the remainder of this section we give a brief introduction to the L\'evy process ideas we will need.

\subsection{L\'evy Process Setup}\label{ssLPS}
We consider a real valued L\'evy process $(X_t)_{t\ge 0}$ 
on a filtered probability space $(\Omega, ({\cal F}_t)_{t\ge 0}, \PP)$,  with canonical triplet $(\gamma,\sigma^2,\Pi)$;  thus, 
 having characteristic function $\EE e^{\rmi \theta X_t}= e^{t\Psi(\theta)}$, $t\ge 0$,
$\theta \in \R$, with exponent
\be\label{ce}
\Psi(\theta)
:= \rmi \theta \gamma - \tfrac 1 2 \sigma^2 \theta^2 
+\int_{\R\setminus\{0\}}\big(e^{\rmi \theta x}-1
-\rmi \theta x{\bf 1}_{\{|x|\le 1\}}\big)\Pi(\rmd x).
\ee
Here $\gamma\in\R$, $\sigma^2\ge 0$ and $\Pi$ is a L\'evy measure on $\R$, i.e., a Borel measure on $\R$ with $\int_{\R\setminus\{0\}}(x^2\wedge 1)\Pi(\rmd x)<\infty$.
The positive, negative and two-sided tails of $\Pi$ are
defined for $x>0$ by 
\be\label{taildef}
\pibar^+(x):=\Pi\{(x,\infty)\},\
 \pibar^-(x):=\Pi\{(-\infty,-x)\},\
 {\rm and}\
 \pibar(x):=\pibar^+(x)+\pibar^-(x).
\ee
Let $\pibarpinv$ denote the inverse function of $\pibar^+$, defined by 
\be\label{defpinv}
\pibarpinv(x)=\inf\{y>0:\pibar^+(y)\le x\}, \ x>0, 
\ee
and  similarly for  $\pibarinv$.
%%,  the inverse function of $\pibar^+$.
%Jumps of $X$ are denoted by $\Delta X_t=X_t-X_{t-}$, $t>0$, with $\Delta X_0=0$,
Throughout, let $\N:=\{1,2,\ldots\}$ and $\N_0:=\{0,1,2,\ldots\}$. 

 Write $(\Delta X_t:=X_t-X_{t-})_{t>0}$, with   $\Delta X_0=0$, for the jump process of $X$, and
 $\Delta X_t^{(1)} \ge \Delta X_t^{(2)}\ge \cdots$
 for the jumps ordered by their magnitudes at time $t>0$.
 Assume throughout that $\Pi\{(0,\infty)\}=\infty$,
%% $\Pi\{\R\setminus\{0\}\}=\infty$, 
 so there are infinitely many positive jumps, a.s., in any 
right neighbourhood of 0.
%% Further assume $\pibar^+(0+)>0$, so 
Thus the $\Delta X_t^{(i)}$ are positive a.s. for all $t>0$ but 
$\lim_{t\dto 0} \Delta X_t^{(i)}=0$ for all $i\in\N$.
 Our objective %%in Section \ref{s5} 
 is to study the ``one-sided trimmed process", by which we mean $X_t$ minus its large positive jumps, at a given time $t$. Thus, the one-sided $r$-trimmed version of  $X_t$ is
\begin{equation}\label{trims}
{}^{(r)}X_t:= X_t-\sum_{i=1}^r {\Delta X}_t^{(i)}, \ r\in\N, \ t>0
\end{equation}
(and we set $^{(0)}X_t\equiv X_t$).
  Detailed definitions and properties of this kind of ordering  and trimming are given in \cite*{bfm2016},
%%Sections \ref{sec: PPP} and \ref{s3}, 
where we identify the positive $\Delta X_t$ with the points of a Poisson point process %%$(\Delta_t)_{t>0}$ 
on $[0,\infty)$. 

Our main result, in Theorem \ref{stable}, is to show that ratios formed by dividing  ${}^{(r)}X_t$, possibly after centering,  by its ordered positive jumps, converge to the corresponding stable ratios when $X$ is in the domain of attraction of a non-normal stable law.
%%In this way we generalise and extend a number of known results.

\section{Convergence of L\'evy Ratios to Stable Limits}\label{s5}
Throughout,   $X$ will be assumed to be 
%%driftless subordinator  
in the domain of attraction of a non-normal stable  random variable at 0 
(or at $\infty$).\footnote{The convergences in this section can be worked out as $t\dto 0$ or as $t\to\infty$. For definiteness and in keeping with modern trends in the area we supply the versions for $t\dto 0$, but little modification is needed for the case $t\to\infty$.}
By this we mean that there are nonstochastic functions $a_t\in\R$ and $b_t>0$
such that $(X_t-a_t)/b_t\todr Y$, for an a.s.  finite random variable $Y$, not degenerate at a constant, and not normally distributed, as $t\downarrow0$. % ($t\to \infty$); equivalently,
The L\'evy tail $\pibar(x)$ is then regularly varying of index $-\alpha$ at 0,  %%(or at $\infty$), 
and the balance conditions
\be\label{bal}
\lim_{x\dto 0} \frac{\pibar^\pm(x)}{\pibar(x)}= a_\pm,
\ee
where $a_++a_-=1$, are satisfied.
If this is the case then the limit  random variable  $Y$ must be a stable  random variable  with index $\alpha$ in $(0,2)$.
We  consider one-sided (positive) trimming, so we always assume $a_+>0$,
and then also $\pibar^+(x)$ is   regularly varying   at 0 with index 
$-\alpha$, $\alpha\in(0,2)$.

Denote by $RV_0(\alpha)$ ($RV_\infty(\alpha)$) the regularly varying functions of index $\alpha\in\R$ at 0 (or $\infty$).
%We know that $\pibar(x)\in RV_{0/\infty}(\alpha)$ iff
%  $\pibarinv(x)\in RV_{\infty/0}(1/\alpha)$.
  %%, where $1/0\equiv \infty$.
%  $RV_{0/\infty}(0)$ are the slowly varying functions at $0$ or $\infty$, and $RV_{0/\infty}(\infty)$ are the rapidly varying functions at $0$ or $\infty$.
When $\pibar^+(\cdot) \in RV_0(-\alpha)$ with $0<\alpha<\infty$ or, equivalently,  the inverse function $\pibarpinv(\cdot) \in RV_\infty(-1/\alpha)$
 (e.g. Bingham, Goldie and Teugels \citeyearpar[Sect. 7, pp.28-29]{BGT87}),  we have the easily verified convergence
\be\label{1b}
t\pibar^+(u\pibarpinv(y/t))
\sim
\frac{\pibar^+(uy^{-1/\alpha}\pibarpinv(1/t))} 
{\pibar^+(\pibarpinv(1/t))}
 \to u^{-\alpha} y \ \text{as} \ t \dto 0,\ {\rm for\ all}\ u, y > 0.
\ee
%When $\alpha=\infty$ we interpret $u^{-\alpha}$  as 
%$\infty.{\bf 1}_{\{0<u<1\}}+1.{\bf 1}_{\{u=1\}}
%+0.{\bf 1}_{\{u>1\}}$.
For $r>0$ write
\ben
\PP(\Gamma_r\in\rmd x)= \frac{x^{r-1} e^{-x}\rmd x}{\Gamma(r)} {\bf 1}_{\{x>0\}},
\een
for the density of $\Gamma_r$,  a Gamma$(r,1)$ random variable, 
which should not be confused with the Gamma function,
$\Gamma(r)=\int_0^\infty x^{r-1}e^{-x}\rmd x$.
 Denote the Beta random variable on $(0,1)$  with parameters $a, b>0$ by ${\rm B}_{a, b}$, 
having density function
\[ f_B(x) = \frac{\Gamma(a+b)}{\Gamma(a)\Gamma(b)} x^{a-1}(1-x)^{b-1} = \frac{1}{B(a,b)}x^{a-1}(1-x)^{b-1},\  0<x<1.
\]

Denote by $(S_t)_{t\ge 0}$  a stable process of index $\alpha\in(0,2)$ having  
 %%triplet $(\gamma_S,0,\Lambda)$, with $\gamma_S\in\R$ and 
 L\'evy measure
\be\label{LM}
\Lambda(\rmd x)= \Lambda_S(\rmd x)=
-\rmd(x^{-\alpha}) {\bf 1}_{\{x>0\}}+ (a_-/a_+)\rmd (-x)^{-\alpha}{\bf 1}_{\{x<0\}},\ x\in\R,
\ee
with characteristic 
exponent
%\footnote{It's convenient for us to write the stable characteristic exponent in this way. It can be put in a more usual form by adjusting the centering function $a_t$.}
\be\label{ces}
\Psi_S(\theta):=\int_{\R\setminus\{0\}}\left(e^{\rmi \theta x}-1
-\rmi \theta x{\bf 1}_{\{|x|\le 1\}}\right)\Lambda(\rmd x),
\ee
and by  $(\Delta S_t:=S_t-S_{t-})_{t>0}$ the jump process of $S$.
Let 
\ben
\Delta S_t^{(1)}\ge \Delta S_t^{(2)}\ge\cdots\ge \Delta S_t^{(n)}\ge \cdots
\een
be the ordered stable jumps at time $t>0$. These are uniquely defined a.s. (no tied values a.s.) since the L\'evy measure of $S$ has no atoms.
The positive and negative tails of $\Lambda$ 
are $\lambar^+(x):=\Lambda\{(x,\infty)\}=x^{-\alpha}$ and $\lambar^-(x):=\Lambda\{(-\infty, -x)\}=(a_-/a_+)x^{-\alpha}$, for $x>0$.  Since $\lambar^+(0+)=\infty$,  the $\Delta S_t^{(i)}$ are positive a.s., $i=1,2,\ldots$, but tend to 0 a.s. as $t\dto 0$.

%L\'evy process in the domain of attraction of a stable  random variable, 
%%%$S=S_1$,  
%having triplet
%$(\gamma,\sigma^2,\Pi)=(\gamma,0,\Pi)$, where the L\'evy measure $\Pi$  attributes no mass to $(-\infty,0)$.
%
%%% This means $X\in bv$ and 
% We additionally assume $X$ has drift $\rmd_X=0$ hence in fact it is a driftless subordinator.

% Since there is no negative component to $\Pi$ we simply write 
% $\pibar(x)$ and $\pibarinv(x)$  for $\pibar^+(x)$ and 
% $\pibarpinv(x)$ from now on.
% The trimmed stable  random variable  ${}^{(r)}S_1$ and the stable process jumps are as in the previous section.

Define a centering function $\rho_X(\cdot)$ for $X$  by 
\bea\label{defnewu}
\rho_X(w):=
\begin{cases}
\displaystyle{
\gamma-\int_{[-1,-w)\cup [w,1]}x\Pi(\rmd x),
%%\gamma-\int_{[w,1]}x\Pi(\rmd x),
}
&0< w\le 1,\\
\\
\gamma+\int_{[-w,-1)\cup (1,w)}x\Pi(\rmd x),
&w>1,
\end{cases}
\eea
and similarly for $\rho_S(w)$, but  with $\gamma$ taken as 0 and $\Lambda$ replacing $\Pi$ in that case.
%and let
%$\wt\rho_t^{(X,r)}:=  \rho(\Delta X_t^{(r)})$ and 
%$\wt\rho_1^{(S,r)}:= \rho(\Delta S_1^{(r)})$, $r\in\N$
%(the tilde is to emphasise that these are random objects).

%The main theorem in the paper, 
%Theorem \ref{stable}, 
%proves joint convergence of the trimmed process $(^{(r)}X_t)$, normalised by those large jumps trimmed,  and after centering, to corresponding stable ratios.
 To state Theorem \ref{stable}, we need some further notation. 
For each $n=2,3,\ldots$ and $0<u<1$, suppose random variables 
$J_{n-1}^{(1)}(u)\ge J_{n-1}^{(2)}(u)\ge\cdots\ge J_{n-1}^{(n-1)}(u)$ are distributed like the decreasing order statistics of $n-1$ independent and identically distributed (i.i.d.)  random variables $(J_i(u))_{1\le i\le n-1}$, each having the distribution 
\be\label{Jdis0}
\PP(J_1(u)\in\rmd x)= \frac{\Lambda(\rmd x) 
{\bf 1}_{\{1\le x\le 1/u\}}}{1-u^{\alpha}},\ x>0.
\ee
Also let $L_{n-1}^{(1)}\ge L_{n-1}^{(2)}\ge\cdots\ge L_{n-1}^{(n-1)}$ be distributed like the decreasing order statistics of $n-1$ i.i.d. random variables $(L_i)_{1\le i\le n-1}$, each having the distribution 
\be\label{Ldis}
\PP(L_1\in \rmd x)= 
\Lambda(\rmd x) {\bf 1}_{\{x>1\}}.   
\ee

Define
\be\label{3a01}
\psi(\theta) =\int_{(-\infty,1)} \left(e^{\rmi\theta x}-1-\rmi \theta x{\bf 1}_{\{|x|\le 1\}}\right)  \Lambda(\rmd x),
\ \theta\in\R,
\ee
and choose  $\theta_0>0$ such that
%%\footnote{$\psi(0)=0$, so $|\psi(\theta)|<1$ for $|\theta|\le$ some $ \theta_0>0$.} 
$|\psi(\theta)|<1$ for $|\theta|\le \theta_0$
(as is possible since $\psi(0)=0$).
Also 
%recall $\psi(\theta)$ in \eqref{3a01}, and
 define $\phi(\theta, u)=\EE e^{\rmi\theta J_1(u)}$, $\theta\in\R$, with $J_1(u)$ having the distribution in \eqref{Jdis0}:
\be\label{phiJ}
\phi(\theta, u)=(1-u^{\alpha})^{-1}
\int_1^{1/u}e^{\rmi\theta x} \Lambda(\rmd x),\ 0<u<1.
\ee
 Let $W=(W_v)_{v\ge 0}$ be a  L\'evy process on $\R$ with triplet $(0,0,\Lambda(\rmd x){\bf 1}_{(-\infty,1)}$), and
$\Gamma_{r+n}$ a Gamma ($r+n$,1) random variable independent of $W$.

When $n=2,3,\ldots$,   $x_k>0$, $1\le k\le n-1$,  $x_{n}=1$, and  $\theta_k\in\R$, $1\le k\le n$,
 write for shorthand
 \be\label{thep}
x_{n+}=
\sum_{k=1}^nx_{k}\quad {\rm and}\quad
\wt\theta_{n+}=\wt\theta_{n+}(x_1,\ldots,x_n):=  \sum_{k=1}^n \frac{\theta_k}{x_{k}},
\ee
and let $\int_{\bfx^\uparrow\ge 1}$ denote integration over the region $\{x_1\ge x_2\ge \cdots \ge x_{n-1}\ge 1\}\subseteq \R^{n-1}$.
Recall that ${}^{(0)}X\equiv X$.

\begin{theorem}\label{stable}%%[Joint Convergence]
Assume $\pibar\in RV_0(-\alpha)$ for some $0<\alpha<2$ and \eqref{bal}. 

{\rm (i)}\
Then   for each $r\in\N_0$, $n\in\N$, as $t\dto 0$,  we have the joint convergence 
\bea\label{liPD}
&&
\bigg( \frac{{}^{(r)}X_t-t\rho_X(\Delta X_t^{(r+n)})}{\Delta X_t^{(r+1)}}, \ldots, \frac{{}^{(r)}X_t-t\rho_X(\Delta X_t^{(r+n)})}{\Delta X_t^{(r+n)}} \bigg) \cr
&&
\todr
\bigg(\frac{^{(r)}S_1-\rho_S(\Delta S_1^{(r+n)})}{\Delta S_1^{(r+1)}}, \ldots, \frac {^{(r)}S_1-\rho_S(\Delta S_1^{(r+n)})} {\Delta S_1^{(r+n)}}\bigg).  %, \ {\rm as}\ t\dto 0.
\eea
%% for each $r\in\N_0$, $n\in\N$.

{\rm (ii)}\ 
%%%%\yfan{problem with $r = 0$, to be fixed.}
When $r\in\N$, $n=2,3,\ldots$, 
the random vector on the RHS of \eqref{liPD} has characteristic function which can be represented,
for  $\theta_k\in\R$,  $1\le k\le n$,  as
\bea\label{38}
%&&
%\lim_{t\dto 0}
%\EE\bigg(\exp\bigg(\rmi \sum_{k=1}^n 
%\frac{\theta_k ({}^{(r)}X_t-t\rho_X(\Delta X_t^{(r+n)}))}{\Delta X_t^{(r+k)}}\bigg)\bigg)\cr
%&&
%=
&&
\EE\exp\bigg(\rmi \sum_{k=1}^n 
\frac{\theta_k\big({}^{(r)}S_1-\rho_S(\Delta S_1^{(r+n)})\big)
}{\Delta S_1^{(r+k)}}\bigg)=
\cr
&&\cr
&&
\int_{\bfx^\uparrow\ge 1}e^{\rmi\wt\theta_{n+} x_{n+}}
\EE\big(e^{\rmi\wt\theta_{n+} W_{\Gamma_{r+n}}}\big) 
\PP\big(J_{n-1}^{(k)}({\rm B}^{1/\alpha}_{r,n})\in\rmd x_k, 1\le k\le n-1\big),
\eea
%where $\int_{\bfx^\uparrow\ge 1}$ denotes  over the region $\{x_1\ge x_2\ge \cdots \ge x_{n-1}\ge 1\}\subseteq \R^{n-1}$.
where $B_{r,n}$ is  a Beta$(r,n)$ random variable independent of the $(J_i(u))$.
Alternatively, recalling \eqref{3a01}, 
when $\max_{1\le k\le n}|\theta_k|\le \theta_0$ 
%%where $\theta_0$ is such that $|\psi(\theta)|<1$ for $|\theta|\le \theta_0$, 
the RHS of \eqref{38} 
 can be written as
\be\label{39}
\int_{\bfx^\uparrow\ge 1}
%%\int_{(x_1,\ldots, x_{n-1})\in\R^{n-1}}
\frac{e^{\rmi\wt\theta_{n+}x_{n+}}}
{\big(1-\psi(\wt\theta_{n+})\big)^{r+n}}
  \PP\big(J_{n-1}^{(k)}({\rm B}^{1/\alpha}_{r,n})\in\rmd x_k, 1\le k\le n-1\big).
  \ee
  
  When $r=0$, \eqref{38} and \eqref{39} remain true as stated if the rvs  $J_{n-1}^{(k)}({\rm B}^{1/\alpha}_{r,n})$ are replaced respectively by  $L_{n-1}^{(k)}$,
being the order statistics associated with the distribution in \eqref{Ldis}.

{\rm (iii)}\ 
When $r\in\N$, $n\in\N$ we have
%, \eqref{liPD} and \eqref{38} take the form
\bea\label{liPD0}
&&
 \frac{{}^{(r)}X_t-t\rho_X(\Delta X_t^{(r+n)})}{\Delta X_t^{(r+n)}}  \todr
\frac{^{(r)}S_1-\rho_S(\Delta S_1^{(r+n)})}{\Delta S_1^{(r+n)}}, \ {\rm as}\ t\dto 0,
\eea
where, recalling \eqref{phiJ}, 
the random variable on the RHS of \eqref{liPD0} has characteristic function
\be\label{36}
%\lim_{t\dto 0}
%\EE\bigg(\exp
%\bigg(\rmi \theta 
%\frac{{}^{(r)}X_t-t\rho_X(\Delta X_t^{(r+n)})}{\Delta X_t^{(r+n)}}\bigg)\bigg)
%=
%\EE\exp\bigg(\rmi \theta 
%\frac{{}^{(r)}S_1-\rho_S(\Delta S_1^{(r+n)})}{\Delta S_1^{(r+n)}}\bigg)=
\frac{e^{\rmi\theta }}{\left(1-\psi(\theta)\right)^{r+n}}
   \,  \EE\big(\phi^{n-1}(\theta, {\rm B}^{1/\alpha}_{r,n})\big), \ \theta\in\R, \,|\theta|\le \theta_0.
\ee
  When $r=0$, \eqref{liPD0} remains true, as does \eqref{36}, if 
  $\phi(\theta, {\rm B}^{1/\alpha}_{r,n})$  in \eqref{36} 
  is replaced by  $\phi(\theta, 0):=
  \int_1^\infty e^{\rmi\theta x} \Lambda(\rmd x)$.
%%\bea\label{26a}
%\lim_{t\dto 0}
%\EE\big(e^{\rmi \theta({}^{(r)}X_t-t\rho_X(\Delta X_t^{(r)}))/
%\Delta X_t^{(r+1)}}\big)
%&=&
%\EE\big(e^{\rmi \theta({}^{(r)}S_1-\rho_S(\Delta S_1^{(r)}))/
%\Delta S_1^{(r+1)}}\big)
%&=&
%e^{\rmi\theta}\EE\big(e^{\rmi\theta W_{\Gamma_{r+1}}}\big),\,
%\theta\in\R.
%\eea
\end{theorem}

Setting $n=1$ in \eqref{liPD0}, and 
(since ${}^{(r)}X_t/\Delta X_t^{(r+1)}=1+{}^{(r+1)}X_t/\Delta X_t^{(r+1)}$)
replacing $r+1$ by $r$ gives 

\begin{corollary}\label{prx}
%%With the assumptions on $X$ and $S$ so far specified, we have 
For each $r\in\N$,  $\theta\in\R$, $|\theta|\le \theta_0$, 
\bea\label{liPD1}
&&
\frac{{}^{(r)}X_t-t\rho_X(\Delta X_t^{(r)})}{\Delta X_t^{(r)}} \todr
\frac{^{(r)}S_1-\rho_S(\Delta S_1^{(r)})}{\Delta S_1^{(r)}}, \ {\rm as}\ t\dto 0,
\eea
where
\be\label{31}
\EE\big(e^{\rmi \theta ({}^{(r)}S_1-\rho_S(\Delta S_1^{(r)}))/
\Delta S_1^{(r)}}\big)
=\EE\big(e^{\rmi\theta W_{\Gamma_{r}}}\big)
=\frac{1}{\big(1-\psi(\theta)\big)^{r}}.
%%, \ r\in\N, \ |\theta|\le \theta_0.
\ee
Further, %%corollary of \eqref{31} is that
$\big({}^{(r)}S_1-\rho_S(\Delta S_1^{(r)})\big)/\Delta S_1^{(r)}\eqdr
W_{\Gamma_{r}}$, being a Gamma-subordinated L\'evy process, is infinitely divisible for each $r\in\N$. 
\end{corollary}
% 
% Similarly, the RHS of \eqref{26} equals
% $e^{\rmi\theta}\EE\left(e^{\Gamma_{r}(1-\psi(\theta)}\right)$
% in this case.

%%\medskip\noindent {\bf Remark.}\ 
The unwieldy centering functions $\rho_X$ and $\rho_S$ in \eqref{liPD}--\eqref{31} can be simplified in many cases.
%%\footnote{ As usual, $\alpha=1$ is an exceptional case.}
%are easily developed using the regular variation of $\pibar ^\pm$. 
Especially, when $X$ is  a subordinator with drift $\rmd_X$, $\rho_X$ can be replaced by $\rmd_X$, and without loss of generality we can assume $\rmd_X=0$. 
The convergences in \eqref{liPD}--\eqref{liPD1} can then be written in terms of Laplace transforms.
%%We get  \eqref{liPD3} in the next theorem, as well as a converse:
% and The regular variation assumption forces $0<\alpha<1$, and
This case recovers a result proved in  Theorem 1.1 of \cite{KeveiMason2014}:
{\it assume $X$ is a driftless subordinator  in the domain of attraction  (at 0) of a stable random variable with index $\alpha\in(0,1)$. Then  for $r\in\N$ }
\bea\label{liPD3}
&&
\frac{{}^{(r)}X_t}{\Delta X_t^{(r)}} \todr {}^{(r)}Y,\ {\rm as}\ t\dto 0,
\eea
where $^{(r)}Y$  is an a.s. finite non-degenerate  random variable.
From Theorem \ref{stable} we can identify ${}^{(r)}Y$ as having the distribution of $^{(r)}S_1/\Delta S_1^{(r)}$, in our notation.
Kevei and Mason show, conversely,  in this subordinator case, that when \eqref{liPD3} holds with $^{(r)}Y$  a finite non-degenerate  random variable, then %$\pibar(x)\in RV_0(-\alpha)$ for some $\alpha\in(0,1)$, thus, 
$X$ is  in the domain of attraction  (at 0) of a stable random variable with index $\alpha\in(0,1)$.
They also give in their Theorem 1.1 a formula for the Laplace transform of  ${}^{(r)}Y$.
We can  state an equivalent version as: {\it suppose \eqref{liPD3} holds.
Then  \eqref{31} becomes
\be\label{31L}
\EE\big(e^{-\lambda{}^{(r)}S_1/\Delta S_1^{(r)}}\big)
=\EE\big(e^{-\lambda W_{\Gamma_{r}}}\big)
=\frac{1}{\left(1+\Psi(\lambda)\right)^{r}}, \ r\in\N,
\ee
where  now $W=(W_v)_{v\ge 0}$ is a  driftless subordinator with  measure 
$\Lambda(\rmd x){\bf 1}_{(0,1)}$, and }
\ben\label{psdef}
\Psi(\lambda) =
\int_{(0,1)} \big(1-e^{-\lambda x}\big)  \Lambda(\rmd x),\
\lambda>0.
\een

\begin{remark}[Negative Binomial Point Process]\label{rem1}\
{\rm
The form of the Laplace transform in \eqref{31L} suggests a connection with the negative binomial point  process of \cite{Gregoire1984}.
That connection is developed  in detail in \cite{NBB2017},  and also forms the basis for a general point measure treatment when $0\le \alpha \le \infty $ in \cite{IMR1}, which contains a converse proof generalising that of \cite{KeveiMason2014}. 
Those results motivate further investigation of the ``large trimming" properties of general L\'evy processes in the spirit of \cite*{BMR2016}. 
%%But we do not explore this further here.
}
\end{remark}

%\begin{remark}
%{\rm 
%As mentioned, Theorem \ref{stable}  is  a kind of multivariate version for  a general L\'evy of 
%Theorem 1.1 of \cite{KeveiMason2014}, who restrict themselves to subordinators (see the remarks relating to our Eq. \eqref{liPD3}).
%They derive  the limit distribution of the ratio of an $r$-trimmed subordinator to its $(r+1)$-st largest jump, proving also  a converse result, as well  as boundary cases (slow and rapid variation of $\pibar$ at 0).
%}
%\end{remark}

\begin{remark}[Modulus Trimming] \label{rem2}\ 
{\rm
Rather than removing large jumps from $X$ as we do in \eqref{trims},
we can remove jumps large in modulus and obtain analogous formulae and results, with appropriate modifications.
The centering function $\rho_X$ in \eqref{defnewu} should then be changed to 
$\gamma-\int_{[-1,-w]\cup [w,1]}x\Pi(\rmd x)$ 
when $0<w\le 1$,   and to
$\gamma+\int_{(-w,-1)\cup (1,w)}x\Pi(\rmd x)$ when $w>1$, 
 and similarly for $\rho_S$.
The norming in Theorem \ref{stable} would then be by jumps large in modulus rather than by large (positive) jumps, and the convergence would be to the analogous modulus trimmed stable process.
The identities in Section \ref{s3} required for the modified proofs can be obtained from analogous formulae for modulus trimming in \cite*{bfm2016}.
}
\end{remark}

\begin{remark}[Connection with $\PD_\alpha^{(r)}$]\label{rem3}\
{\rm
When $X$ is a driftless subordinator, we obtain from \eqref{liPD} with $n\in\N$ that
\be\label{liPD2}
\bigg( \frac{\Delta X_t^{(r+1)}}{{}^{(r)}X_t}, \ldots, \frac{\Delta X_t^{(r+n)}}{{}^{(r)}X_t} \bigg) 
\todr
\bigg(\frac{\Delta S_1^{(r+1)}}{{}^{(r)}S_1}, \ldots, \frac{\Delta S_1^{(r+n)}}{{}^{(r)}S_1}\bigg), \ {\rm as}\ t\dto 0.
\ee
When $n\to \infty$, the $n$-vector on the RHS tends to a vector $(V^{(r)}_1,V^{(r)}_2, \ldots)$ on the infinite simplex with the generalised  Poisson-Dirichlet distribution $\PD_\alpha^{(r)}$ defined in \cite{trimfrenzy}.
When $r = 0$, this reduces to the Poisson-Dirichlet distribution generated from the stable subordinator, denoted by $\PD(\alpha, 0)$  in \cite{PY1997}, which was first noted by \cite{Kingman1975}. 
}
\end{remark}

To complete this section we continue to consider the case when $X$ is a driftless subordinator.
% Next is a converse result to \eqref{liPD3}.  %Theorem \ref{stable}.
Our final result shows that ratios of the form ${}^{(r+n)}X_t/\Delta X_t^{(r)}$  as in \eqref{liPD2} have strong stability properties.
In the next theorem the interesting aspect is the uniformity of convergence in neighbourhoods of 0;
although $\Delta X_t^{(r)}\dto 0$ a.s. as $t\dto 0$, the remainder after removing an increasing number of jumps, $r+n$, from $X$ is shown to be small order $\Delta X_t^{(r)}$, in probability as $n\to\infty$, uniformly on compacts. %% intervals.
%% $t\in(0,t_0]$, $t_0>0$.

\begin{theorem}\label{prin}
Suppose $X$ is a driftless subordinator with $\pibar\equiv\pibar^+\in RV_0(-\alpha)$ for some $0<\alpha<1$. Then  for each $r\in\N$
\be\label{B15}
\frac{{}^{(r+n)}X_t}{\Delta X_t^{(r)}}\topr 0,\ {\rm as}\ n\to\infty,
\ee
uniformly in $t\in(0,t_0]$, for any $t_0>0$.
\end{theorem}

\begin{remark}
{\rm 
By the uniform in probability convergence in Theorem \ref{prin} we mean 
\be\label{uip}
\lim_{n\to \infty}\PP({}^{(r+n)}X_t>\veps \Delta X_t^{(r)})=0,\
{\rm uniformly\ in} \ 0<t\le t_0, \ {\rm for\ all}\ \veps>0.
\ee
Since ${}^{(r+n)}X_t$ is monotone in $n$, this is equivalent to a kind of ``uniform almost sure" convergence, as follows.
With ``i.o." standing for ``infinitely often", and $\veps>0$, $t>0$, 
\bean
\PP({}^{(r+n)}X_t>\veps \Delta X_t^{(r)}\ {\rm i.o.},\ n\to \infty)
&=&
\lim_{m\to\infty}
\PP({}^{(r+n)}X_t>\veps \Delta X_t^{(r)}\ {\rm for\ some}\ n>m)\cr
&\le &
\lim_{m\to\infty}
\PP({}^{(r+m)}X_t>\veps \Delta X_t^{(r)})
=0,
\eean
where the convergence is uniform in $0<t\le t_0$, by \eqref{uip}.
}
\end{remark}
%%\newpage

\section{Representations for Trimmed L\'evy Processes}\label{s3}
% In the next  sections we will identify the $\Delta_t^{(r)}$ with the ordered jumps of a L\'evy process $X=(X_t)_{t\ge 0}$ on $\R$.
%%, with L\'evy measure restricted to $(0,\infty)$, thus, having no negative jumps. 
%We borrow the general setup from BFM \citeyearpar{bfm2016}.
%%Consider  
%Although in Sections \ref{s5} and \ref{ps5} we deal with subordinators, much of our setup and basic results are quite general. 
In the present section we revert to considering an arbitrary  real valued L\'evy process $(X_t)_{t\ge 0}$, set up as in Section \ref{intro}
(see \eqref{ce} and \eqref{taildef}), and derive the identities required for the proofs of the results in Section \ref{s5}.
Fundamental to these identities is a general representation for the joint  distribution of 
${}^{(r)}X_t$ and its large jumps, given in  \cite{bfm2016},
%BFM \citeyearpar{bfm2016},
%%$\Delta X_t^{(r)})$,
which allows for possible tied values in the jumps.\footnote{A different but equivalent distributional representation when $X$ is a subordinator is in Proposition 1 of \cite{KeveiMason2013}.}
%Kevei, P. and  Mason, D.M. (2013)
%Randomly weighted self-normalized L\'evy processes,
%Stochastic Processes and their Applications 123 (2013) 490--522.
Our main theorem in this section, Theorem \ref{P1}, applies it to
 derive formulae for the conditional distributions 
 of the trimmed L\'evy given some of  its large jumps.
We expect these formulae will have useful applications in other areas too.

%We need a formula for the distribution of the trimmed L\'evy together with its large jumps, given in \cite{bfm2016}.
To state the  \cite{bfm2016}  representation, recall the definition of 
%%let $\pibarpinv(x)=\inf\{y>0: \pibar^+(y) \le x\}$,  $ x>0$, be 
the right-continuous inverse  $\pibarpinv(x)$ of $\pibar^+$ in \eqref{defpinv}, and for each $v>0$ introduce a L\'evy process $(X_t^v)_{t\ge 0}$ 
%as 
%\begin{equation}\label{2.2b}
%X_t^v:=X_t-\sum_{0<s\le t} \Delta X_s\;\bfeins_{\{\Delta X_s\ge \pibar^{\leftarrow}(v)\}},\ t>0,
%\end{equation}
%where $\pibarpinv(x)=\inf\{y>0: \pibar^+(y) \le x\}$,  $ x>0$, is the right-continuous inverse of $\pibar^+$. 
%Then $(X_t^v)_{t\ge 0}$ is a well defined L\'evy process 
%%
having the canonical triplet 
\bea\label{trip3}
&&
\left(\gamma^v, \sigma^2,\,  \Pi^v(\rmd x)\right):=\cr
&&
\bigg(\gamma-{\bf 1}_{\{\pibarpinv(v)\le 1\}}\int_{\pibarpinv(v)\le x \le 1}x \Pi(\rmd x),\,  \sigma^2,\,  \Pi(\rmd x){\bf 1}_{\{x<\pibarpinv(v)\}}\bigg).
\eea
Further, let  $G_t^{v}=\pibarpinv(v)Y_{t\kappa(v)}$ for $v>0$, $t>0$, with
 $\kappa(v):=\pibar^+(\pibarpinv(v)-)-v$ and $(Y_t)_{t\ge 0}$ a  homogeneous Poisson process with $\EE(Y_1)=1$, independent of  $(X_t^v)_{t\ge 0}$.
Let  $r\in \N$ and recall that  $(\SSSS_i)$ are Gamma$(i,1)$ random variables, $i\in\N$.
 Assume that  $(X_t^v)_{t\ge 0}$,  $(G_t^{v})_{t\ge 0}$ and $(\SSSS_i)$ are independent as random elements for each $v>0$. 
Then Theorem 2.1, p.2329,  together with  Lemma 1, p.2333, of \cite{bfm2016} give, for each $t>0$,
$r,m\in \N$, $1\le m\le r$, 
% \begin{equation}\label{2rrep_1}
% \left(^{(r)}X_t,\, \Delta X_t^{(r)}\right)
% \eqdr \left(X_t^{v} + G_t^v,\,
% \pibarinv\left(v\right)\right) \bigg|_{v=\SSSS_r/t}.
% \end{equation}
 \be\label{2333}
\big(^{(r)}X_t,\, \Delta X_t^{(m)},\ldots,  \Delta X_t^{(r)}\big)
\eqdr
\big(X_t^{\Gamma_{r}/t}+G_t^{\Gamma_{r}/t},\, 
\pibarpinv(\Gamma_m/t),\ldots,\pibarpinv(\Gamma_{r}/t)\big).
\ee
%%Note that then the expression in \eqref{cor1bb} is zero when $\Pi$ has no atoms in $(0,\infty)$. %%But we do not assume this.
%in \eqref{trip3}, the  Poisson process $Y$, and the quantities $G_t^v$ and $\kappa$ as in the statement of Theorem \ref{P1}. 

 We need some further notions.
For each $y>0$ introduce another L\'evy process $(X_t^{(y)})_{t\ge 0}$ 
having the canonical triplet 
\be\label{tripy}
\left(\gamma^{(y)}, \sigma^2,\,  \Pi^{(y)}(\rmd x)\right):=
\bigg(\gamma-{\bf 1}_{\{y\le 1\}}\int_{y\le x \le 1}x \Pi(\rmd x),\,  \sigma^2,\,  \Pi(\rmd x){\bf 1}_{\{x<y\}}\bigg),
\ee
and another process   $(G_t^{(y,v)})$ defined such that  $G_t^{(y,v)}=yY_{t\kappa(y,v)}$ for $y,v,t>0$, where again $(Y_t)_{t\ge 0}$ is a  homogeneous  Poisson process with $\EE(Y_1)=1$, now independent of  $(X_t^{(y)})_{t\ge 0}$, and  $\kappa(y,v):=\pibar^+(y-)-v$.
 
 We need to distinguish situations when $\pibar^+$ is or is not continuous at a point. 
 %% $y$.   %% $y_1\ge y_2\ge \cdots\ge y_r>0$ 
%% To save writing, l
 Let $A_\Pi$ denote the points of discontinuity of $\pibar^+$ in $(0,\infty)$.
When $y_i\in A_\Pi$, 
%%points of discontinuity of $\pibar^+$, 
set
 \be\label{defab}
 a_i=a_i(y_i)=\pibar^+(y_i) <  b_i=b_i(y_i)=\pibar^+(y_i-).
 \ee
 Note that $\pibarpinv(v)$ takes the same value, namely, $y_i$, 
 for any $v\in [a_i,b_i)$.
%%Use \eqref{defab} and $(G_t^{(y,v)})$ to 
When $r,m\in \N$ with $1\le m\le r$ and $y_r\in A_\Pi$, define conditional expectations
\be\label{Krdef}
K_{m,r}(\theta,t,y_m,\ldots, y_r)=
\EE\Big(e^{ \rmi\theta G_t^{(y_r,\Gamma_r/t)}}\,
\big|\,\Gamma_i/t\in [a_i,b_i),\, m\le i\le r\Big),
\ee
% and 
% \be\label{K1def}
% K_{r,1}(\theta,t,y_r)=
% \EE\Big(e^{ \rmi\theta G_t^{(y_r,\Gamma_r/t)}}\,
% \big|\,\Gamma_r/t\in (a_r,b_r)\Big),
% \ee
for $t>0$, $\theta\in\R$. 
When $y_r\notin A_\Pi$, i.e.,  $\pibar^+(y_r)=\pibar^+(y_r-)$, 
 we set $ G_t^{(y_r,\cdot)}=0$ and in this case we understand
$ K_{m,r}(\theta,t,y_m,\ldots, y_r)=1$.  %%=K_{r,1}(\theta,t,y_r)$.
 When  $y_r\in A_\Pi$ but $y_i\notin A_\Pi$ for one or more $i$, $m\le i<r$,  
 we understand the corresponding events $\{\Gamma_i/t\in [a_i,b_i)\}$ are omitted from the conditioning in \eqref{Krdef}. 
 
With this notation in place we can now state  Theorem \ref{P1}, the main result of this section, which 
provides  in characteristic function  form 
the conditional  distribution of
$^{(r)}X_t$, given $\Delta X_t^{(r)}$, or  given  $\Delta X_t^{(m)},\ldots,  \Delta X_t^{(r)}$.

\begin{theorem} \label{P1}  
%%%%Assume $\pibar^+(0+)=\infty$.
Take  integers $r,m\in \N$ with $1\le m\le r$, and real numbers
 $y_m\ge \cdots\ge y_r>0$, $\theta\in\R$,  $t>0$. Then we have the identities
\be\label{q4}
\EE\big(e^{ \rmi\theta {}^{(r)}X_t}\,\big|\,\Delta X_t^{(r)}=y_r\big) =
\EE\big(e^{ \rmi\theta X_t^{(y_r)}}\,\big)
K_{r,r}(\theta,t,y_r)
\ee
and 
\be\label{q0}
\EE\big(e^{ \rmi\theta {}^{(r)}X_t}\,\big|\,\Delta X_t^{(i)}=y_i,\, m\le i\le r\big) =
\EE\big(e^{ \rmi\theta X_t^{(y_r)}}\,\big)
K_{m,r}(\theta,t,y_m,\ldots,y_r).
\ee
%%s at points of increase $v>0$ of $\pibar^+$.
%In \eqref{q4},  $(X_t^v)_{t\ge 0}$ is a L\'evy process, indexed by $v>0$,  
%%as 
%%\begin{equation}\label{2.2b}
%%X_t^v:=X_t-\sum_{0<s\le t} \Delta X_s\;\bfeins_{\{\Delta X_s\ge \pibar^{\leftarrow}(v)\}},\ t>0,
%%\end{equation}
%%where $\pibarpinv(x)=\inf\{y>0: \pibar^+(y) \le x\}$,  $ x>0$, is the right-continuous inverse of $\pibar^+$. 
%%Then $(X_t^v)_{t\ge 0}$ is a well defined L\'evy process 
%%%
%having canonical triplet
\end{theorem}
%
%%%\medskip\noindent {\bf Remark}\ 
%\begin{remark}\label{rem5}{\rm 
%The RHS of \eqref{q4} does not depend on $r$, so it also equals
%$\PP\big({}^{(1)} X_t\le x\,\big|\, \Delta X_t^{(1)}=\pibarpinv(v)\big)$.
%Similarly in \eqref{B9} and \eqref{B10} below.
%}\end{remark}

\medskip\noindent {\bf Proof of  Theorem \ref{P1}:}\
We prove \eqref{q4}, then show how it can be extended to \eqref{q0}.
First suppose $y_r\in A_\Pi$.  %% is an atom of $\Pi$.
From \eqref{2333} we have
\bea\label{Y1}
&&\PP(\Delta X_t^{(r)} = y_r) =
\PP\big(\pibarpinv(\Gamma_r/t)=y_r\big)\cr
&=&
\PP(\Gamma_r/t \in [\pibar^+(y_r), \pibar^+(y_r-))) 
=
\PP(\Gamma_r/t \in [a_r,b_r))> 0.
\eea
(In the last equality, recall \eqref{defab}.)
Since the probability in \eqref{Y1} is positive, we can compute, by elementary means, 
using  \eqref{2333} again,
\bea\label{Y2}
&&\PP({}^{(r)}X_t \le x \, \big|\,  \Delta X_t^{(r)} = y_r) 
= \frac{\PP({}^{(r)}X_t \le x, \,  \Delta X_t^{(r)} = y_r)}{\PP(\Delta X_t^{(r)} = y_r)} \cr
% &&=
% \frac{\PP \Big( X_t^{\Gamma_r/t}  + G_t^{\Gamma_r /t} \le x,   \pibarpinv(\Gamma_r/t) = y_r\Big)}{\PP(\Gamma_r/t \in [a_r,b_r))} \cr
&&=\frac{\PP \Big( X_t^{\Gamma_r/t}  + G_t^{\Gamma_r /t} \le x,   \Gamma_r/t \in [a_r,b_r)\Big)}{\PP(\Gamma_r/t \in [a_r,b_r))}\cr
&&=
\PP \Big( X_t^{\Gamma_r/t}  + G_t^{\Gamma_r /t} \le x\big|
\Gamma_r/t \in R(y_r)\Big), 
\eea
where $R(y_r):= [a_r, b_r)$.
%%%which is exactly (3.2') suggested by the referee.
Going over to characteristic functions, we find, since 
  $X_t^v$,  $G_t^{v}$ and $\SSSS_r$ are independent for each $v>0$,
\be\label{Y3}
\EE\big(e^{ \rmi\theta {}^{(r)}X_t}\,\big|\,\Delta X_t^{(r)}=y_r\big) 
=
\int_{v\in R(y_r)}
\EE(e^{\rmi\theta X_t^{v}}) \, \EE(e^{\rmi\theta G_t^{v}})\frac{\PP\left(\Gamma_r/t\in \rmd v\right)}{\PP(\Gamma_r/t \in  R(y_r))}.
\ee
 Whenever $v\in R(y_r)$,  then $\pibarpinv(v)=y_r$ and $X_t^v=X_t^{(y_r)}$ (recall \eqref{tripy}), 
 while $\kappa(v)=\pibar^+(y_r-)-v=\kappa(y_r,v)$ and 
 $G_t^v=y_rY_{t\kappa(y_r,v)}= G_t^{(y_r,v)}$. Consequently the RHS of \eqref{Y3} is 
\ben %%\label{Y4}
\EE(e^{\rmi\theta X_t^{(y_r)}})
\EE\big(e^{\rmi\theta G_t^{(y_r,\Gamma_r/t)}}\big|\Gamma_r/t\in R(y_r)\big)
%%{\PP(\Gamma_r/t \in  R(y_r))} 
=
\EE\big(e^{ \rmi\theta X_t^{(y_r)}}\big)
K_{r,r}(\theta,t,y_r),
\een
as required for \eqref{q4}.

The conditional probability in \eqref{Y2} is in fact 
the Radon-Nikodym derivative of the measure 
$\PP\big({}^{(r)} X_t\le x,\, \Delta X_t^{(r)}\le \cdot\big)$
with respect to the measure 
$\PP\big(\Delta X_t^{(r)}\le \cdot\big)$ on $(0,\infty)$
%%which is elementary to calculate 
when  $y_r$ is an atom of $\Pi$.
Alternatively, suppose  $\Pi$ is continuous at $y_r$.
Then we write, from \eqref{2333}, 
for $t>0$, $y_r>0$, 
\be\label{q5}
\PP\big({}^{(r)} X_t\le x,\, \Delta X_t^{(r)}\le y_r\big)
=
\int_{\{v>0:\,\pibarpinv(v)\le y_r\}}
 \PP\left(X_t^{v} + G_t^{v}\le x\right)\PP\left(\Gamma_r/t\in \rmd v\right)
 \ee
 and
 \be\label{q6a}
\PP\big(\Delta X_t^{(r)}\le y_r\big)
=
\int_{\{v>0:\, \pibarpinv(v)\le y_r\}}\PP\left(\Gamma_r/t\in \rmd v\right).
%%=\PP\big(\pibarpinv(\Gamma_r/t)\le y_r\big).
 \ee
 Since $\PP\left(\Gamma_r/t\in \cdot\right)$ is absolutely continuous with respect to Lebesgue measure we can use the differentiation formula in Thm.2, p.156 of \cite{zaanen1958} to calculate the  Radon-Nikodym derivative. 
 Thus we evaluate \eqref{q5} and \eqref{q6a} over intervals
 $(y_r-\veps^-,y_r+\veps^+)$ and take the limit of the ratio as $
 \veps^\pm\dto 0$. 
This produces
\be\label{Y5}
\PP({}^{(r)}X_t \le x \, \big|\,  \Delta X_t^{(r)} = y_r) 
=\PP(X_t^{(y_r)}\le x),
\ee 
 and since $K_{r,r}(\theta,t,y_r)=1$ in this case, we get \eqref{q4} again.

 This completes the proof of \eqref{q4}. Next we extend it to \eqref{q0}.
Assume  $y_m\ge \cdots \ge y_r>0$ are in $A_\Pi$.
Then \eqref{Y1} generalises straightforwardly  to 
\bea\label{Y6}
\PP(\Delta X_t^{(i)} = y_i, m\le i \le r)
=
\PP(\Gamma_i/t \in [a_i,b_i),\, m\le i\le r)> 0,
\eea
and \eqref{Y2} becomes
\bea\label{Y7}
&&\PP({}^{(r)}X_t \le x \, \big|\,  \Delta X_t^{(i)} = y_i,\, m\le i\le r) \cr
&&=
\PP \Big( X_t^{\Gamma_r/t}  + G_t^{\Gamma_r /t} \le x\big|
\Gamma_i/t \in [a_i,b_i),\, m\le i\le r\Big).
\eea
Going over to characteristic functions and recalling  
$K_{m,r}(\theta,t,y_m,\ldots, y_r)$ in \eqref{Krdef} 
we get \eqref{q0}.

The cases  when some or all of the $y_i$ are continuity points of $\Pi$
 can be analysed as for \eqref{q4}.    Since we do not need these formulae for the proofs we omit details.
  \halmos
  
  \begin{remark}
  {\rm (i)\
  When calculating conditional probabilities, we should check that they have the requisite measurability and integrability properties.
  The expressions in \eqref{Y2} and \eqref{Y7} are clearly measurable with respect to their variables,
  and that they integrate to give the respective joint distributions of 
  ${}^{(r)} X_t$ and the relevant $\Delta X_t^{(i)}$ is easily checked by decomposing integrals into discrete and  absolutely continuous components.  
  Effectively, since all of our calculations ultimately involve integration with respect to the absolutely continuous gamma distributions, the needed  properties follow automatically.
  
(ii)\ \eqref{Krdef}, \eqref{q4} and \eqref{q0} show that in general 
the Markov property for the ordered large jumps
does not hold, as  $K_{1,r}(\theta,t,y_1,\ldots, y_r)\ne K_{r,r}(\theta,t, y_r)$ in general.
But when $y_r$ is a continuity point of $\pibar^+$,  then equality does hold here and we do have the Markov property.
This parallels the similar situation for order statistics of i.i.d. random variables.
  }     \end{remark}

Using Theorem \ref{P1} and \eqref{tripy}, conditional characteristic functions of ${}^{(r)}X_t$ can be written as in the next corollary.
For \eqref{B80} and \eqref{B81}, set $m=1\le r$ in \eqref{q0}.

\begin{corollary}\label{conL} 
%%Assume $\pibar^+(0+)=\infty$.
For $r\in\N$, $y_1\ge y_2\ge \cdots\ge y_r>0$, $\theta\in\R$,  $t>0$,
\bea \label{B80}
&&
\EE\big(e^{ \rmi\theta {}^{(r)}X_t}\,
\big|\,\Delta X_t^{(i)}=y_i,\, 1\le i\le r\big)\cr
&&=
\exp\Big(\rmi\theta t\gamma^{(y_r)}- \tfrac{1}{2}t\sigma^2\theta^2
+t\int_{(-\infty,y_r)}\big(e^{\rmi\theta x}-1-\rmi\theta x{\bf 1}_{\{|x|\le 1\}}\big)\Pi(\rmd x)\Big)\cr
&&
\hskip6.5cm \times
K_{1,r}(\theta,t,y_1,y_2,\ldots,y_r).
\eea

Suppose $X$ is  a subordinator  
(so $\sigma^2=0$) with drift
$\rmd_X:=\gamma-\int_{0<x\le 1}x\Pi(\rmd x)$. 
Then the RHS of \eqref{B80} can be replaced by
\be\label{B81}
\exp\bigg(\rmi\theta t\rmd_X
+t\int_{(0,y_r)}\big(e^{\rmi\theta x}-1\big)\Pi(\rmd x)\bigg)
\times K_{1,r}(\theta,t,y_1,y_2,\ldots,y_r).
\ee
\end{corollary}
%
%\begin{remark}\label{rem5}
%{\rm 
%Note that the RHS of \eqref{B80} does not depend on $r$, so it also equals
%$\EE\left(e^{\rmi\theta\, {}^{(1)}X_t}
%\bigg|\Delta X_t^{(1)}=u_r\right)$.
%A similar comment applies to \eqref{B9} and \eqref{B10} below.
%}
%\end{remark}
The next corollary follows immediately from \eqref{q0}.
Recall the definition of $\rho_X$ in \eqref{defnewu}.
For \eqref{B9}, replace $r$ by $r+n$ and set $m= r$ in \eqref{q0}.
%Define
%\bea\label{defnewu}
%\rho(u):=
%\begin{cases}
%\displaystyle{
%\gamma-\int_{[u,1]}x\Pi(\rmd x),
%}
%& u\le 1,\\
%\\
%\gamma+\int_{[-u,-1)\cup (1,u)}x\Pi(\rmd x),
%&u>1.
%\end{cases}
%\eea

\begin{corollary}\label{pr}
For $r\in\N$, $n\in\N_0$, $y_r\ge \cdots\ge y_{r+n}>0$, $\theta\in\R$,  $t>0$,
\begin{align}\label{B9}
&
\EE\bigg(\exp\Big(\rmi\theta
\frac{{}^{(r+n)}X_t-t\rho_X(\Delta X_t^{(r+n)})}{\Delta X_t^{(r+n)}}\Big)
\Big|\Delta X_t^{(k)}=y_k,r\le k\le r+n\bigg)
=
e^{-t\sigma^2\theta^2/2y_{r+n}^2}  \times    \nonumber \\
&
\exp \Big( t\int_{(-\infty,1)}
\big(e^{\rmi\theta x}-1-\rmi\theta x{\bf 1}_{\{|x|\le 1\}}\big)
\Pi\big(y_{r+n}\rmd x\big) \Big)
\times K_{r,r+n}(\theta/y_{r+n},t,y_r,\ldots,y_{r+n}).
\end{align}
Suppose $X$ is a subordinator  with drift $\rmd_X$.
Then   %%\eqref{B9} can be replaced by %%for $r\ge 2$ by
\bea\label{B10}
&&\EE\bigg(\exp\Big(\rmi\theta\frac{{}^{(r)}X_t-t\rmd_X}{\Delta X_t^{(r)}}\Big)
\Big|\,\Delta X_t^{(i)}=y_i,\, 1\le i\le r\big)
\bigg)\nonumber \\
&&
=\exp\bigg(t\int_{(0,1)}
\big(e^{\rmi\theta x}-1\big)\Pi\big(y_r\rmd x\big)\bigg)
\times K_{1,r}(\theta/y_r,t,y_1,y_2,\ldots,y_r).
\eea
\end{corollary}

\medskip\noindent {\bf Proof of  Corollaries \ref{conL} and  \ref{pr}:}\
\eqref{B80} follows from Theorem \ref{P1}, using \eqref{tripy}.
Then \eqref{B81} follows from \eqref{B80}
by rearranging the centering terms.
\eqref{B9} follows from \eqref{B80} and \eqref{defnewu},
and \eqref{B10} follows from \eqref{B9}.
\halmos

\medskip
Another formula follows similarly from \eqref{B81}:
\begin{corollary}\label{pr2}
Suppose $X$ is  a subordinator   with drift $\rmd_X$. 
Then for $\theta\in\R$, $t>0$,  $r\in\N$, $n\in\N$,
\begin{align}\label{B11}
&\EE\Big(\exp\bigg(\rmi\theta\frac{{}^{(r+n)}X_t-t\rmd_X}{\Delta X_t^{(r)}}\bigg)
\bigg|\Delta X_t^{(i)}=y_i,\, r\le i\le r+n\Big)\nonumber \\
&=\exp\Big(t\int_{(0,y_{r+n})}
\big(e^{\rmi\theta x/y_r}-1\big)\Pi(\rmd x) \Big)
\times K_{r,r+n}(\theta/y_r,t,y_r,\ldots,y_{r+n}).
\end{align}
\end{corollary}

For the proofs in Section \ref{ps5} we also need the following result. 

%[NOTE WE ALSO NEED TO MODIFY THAT PAPER \cite{stableJump2} IN LINE WITH THE PRESENT ONE.]

\begin{proposition}\label{large_j}%%[Ratios bigger than 1]
Suppose $\pibar(\cdot) \in RV_0(-\alpha)$ with $\alpha\in(0,2)$, and keep $r \in \N$ and $n =2,3,\ldots$. Take $x_k \ge 1$ for $1\le k\le n-1$.
Then for $x>0$
\bea\label{top-}
&&
\lim_{t\dto 0}
\PP\bigg(\frac{\Delta X_t^{(r+k)} }{\Delta X_t^{(r+n)}}>x_k,\ 1\le k\le n-1 \, 
\Big| \,\Delta X_t^{(r+n)}=x\pibarpinv(1/t)\bigg)\cr
&&=\
\PP\Big(J_{n-1}^{(k)}\big({\rm B}^{1/\alpha}_{r,n}\big)>x_k, 1\le k\le n-1 \Big),
\eea
where $J_{n-1}^{(1)}(u)\ge J_{n-1}^{(2)}(u)\ge\cdots\ge J_{n-1}^{(n-1)}(u)$ are the  order statistics 
 associated with the distribution in \eqref{Jdis0}, %%the $J_i(u)$ are as in \eqref{13b} 
 $B_{r,n}$ is a Beta$(r,n)$  random variable  independent of  
$(J_i(u))_{1\le i\le n-1}$,
and  the limit is taken as $t\dto 0$ through points $x$ such that $x\pibarpinv(1/t)$ is a point of decrease of $\pibar^+$.

When $r=0$, \eqref{top-} remains true if the RHS is replaced by 
\be\label{top-L}
\PP(L_{n-1}^{(k)}>x_k, 1\le k\le n-1)
\ee
 where $L_{n-1}^{(k)}$ are the order statistics associated with the distribution in \eqref{Ldis}.
\end{proposition}

\begin{remark}\label{remN}\
{\rm 
\eqref{top-} and \eqref{top-L} can be stated in a unified fashion if we make the convention that $B_{0,n}\equiv 0$ a.s., put $u=0$ in \eqref{Jdis0}, and identify $(J_i(0))$ with  a sequence $(L_i)$ of
independent and identically distributed random variables each having the  distribution in \eqref{Ldis}.
Similarly for the corresponding statements in Theorem \ref{stable}.
}
\end{remark}

\medskip\noindent {\bf Proof of Proposition \ref{large_j}:}\
This is a variant of the proof of  Theorem \ref{P1}.
Assume $\pibar(\cdot) \in RV_0(-\alpha)$ with $\alpha\in(0,2)$ and choose $r\in\N_0$, $n=2,3,\ldots$, $x_k\ge 1$.
For brevity write $q_t:= \pibarpinv(1/t)$, $t>0$.
First suppose $\pibar^+$ is discontinuous at $xq_t$, $x>0$, so
\ben
\PP(\Delta X_t^{(r+n)}=xq_t)=\PP(\Gamma_{r+n}\in[a_t(x),b_t(x)),
\een
where 
$
 a_t(x):=t\pibar^+(xq_t) <  b_t(x):=t\pibar^+(xq_t-)$,
 and consider the ratio
 \bea\label{rat}
 &&
\frac{\PP\big(\Delta X_t^{(r+k)}>x_k\Delta X_t^{(r+n)},\, 1\le k\le n-1,\, \Delta X_t^{(r+n)}=xq_t\big)}
{\PP\big(\Delta X_t^{(r+n)}=xq_t\big)}\cr
&&\cr
&&=
%\frac{\PP\big(\Delta X_t^{(r+k)}>x_kxq_t,\, 1\le k\le n-1,\, a_t(x)\le \Gamma_{r+n} <b_t(x)\big)}
%{\PP\big(a_t(x)\le \Gamma_{r+n}<b_t(x)\big)}.
\frac{\PP\big(\overline \Pi^{+, \leftarrow}(\Gamma_{r+k}/t)>x_kxq_t,\, 1\le k\le n-1,\, a_t(x)\le \Gamma_{r+n} <b_t(x)\big)}
{\PP\big(a_t(x)\le \Gamma_{r+n}<b_t(x)\big)}.\
\eea
 With $f_{r+n}(y)$ as the bounded, continuous, density of $\Gamma_{r+n}$, the denominator in \eqref{rat} is, by the mean value theorem, 
 \be\label{Nu}
 \int_{a_t(x)}^{b_t(x)} f_{r+n}(y)\rmd y
 = (b_t(x)-a_t(x)) f_{r+n}(\xi_t(x)),
 \ee
for some  $\xi_t(x)\in [a_t(x),b_t(x))$.
 Let $c_t(x_k,x):= t\pibar^+(x_kxq_t)$. 
 Recalling \eqref{2333}, the numerator in \eqref{rat} can be written as 
 \bea\label{x1}
 &&
  \PP\big(\Gamma_{r+k}< t\pibar^+(x_kxq_t),\, 1\le k\le n-1,\, a_t(x)\le \Gamma_{r+n} <b_t(x)\big)\cr
  &&=
  \int_{a_t(x)}^{b_t(x)} 
  \PP\big(\Gamma_{r+k}< c_t(x_k,x),\, 1\le k\le n-1\,\big|\Gamma_{r+n}=y\big)
  f_{r+n}(y)\rmd y\cr
  &&\cr
  &=&
   \int_{a_t(x)}^{b_t(x)} 
   \PP\bigg(\frac{\Gamma_{r+k}}{\Gamma_{r+n}}< \frac{c_t(x_k,x)}{y},\, 1\le k\le n-1\bigg)     f_{r+n}(y)\rmd y.
   \eea
 In the last equation we used that
$(\Gamma_{r+k}/\Gamma_{r+n})_{1\le k\le n-1}$ is independent of $\Gamma_{r+n}$ (using ``beta-gamma  algebra"; see, e.g., \citet[p.11]{Pitman2006}).

 Again using the mean value theorem  the last expression in \eqref{x1} equals
 \be\label{DE}
(b_t(x)-a_t(x)) f_{r+n}(\eta_t(x))
   \PP\bigg(\frac{\Gamma_{r+k}}{\Gamma_{r+n}}< \frac{c_t(x_k,x)}{\eta_t(x)},\, 1\le k\le n-1\bigg) 
 \ee
 for some $\eta_t(x)\in [a_t(x),b_t(x))$.
Recall \eqref{1b} and that  $q_t:= \pibarpinv(1/t)$ to see that each of $a_t(x)$, $b_t(x)$, $\xi_t(x)$ and $\eta_t(x)$ tends to $x^{-\alpha}$,  that $f_{r+n}(\xi_t(x))$ and  $f_{r+n}(\eta_t(x))$
both tend to  $f_{r+n}(x^{-\alpha})$, and that $c_t(x_k,x)$ tends to $(x_kx)^{-\alpha}$, all as $t\dto 0$. 
 Take the ratio of the numerator of \eqref{rat} in the form \eqref{DE} to
  the denominator in the form \eqref{Nu}, and let $t\dto 0$
  to get the limit of the ratio in \eqref{rat} as
  \be\label{de}
   \PP\bigg(\frac{\Gamma_{r+k}}{\Gamma_{r+n}}< x_k^{-\alpha},\, 1\le k\le n-1\bigg).
 \ee
 
 This gives an expression for the limits on the lefthand side of \eqref{top-} and \eqref{top-L}. To write them in the forms of the  righthand sides of \eqref{top-} and \eqref{top-L}, first take $r\in\N$, and use the fact that, 
 conditionally on $\Gamma_r/\Gamma_{r+n} =s>0$,  
\be\label{df}
\left(\frac{\Gamma_{r+1}}{\Gamma_{r+n}},\dots, \frac{\Gamma_{r+n-1}}{\Gamma_{r+n}}\right)
\eqdr \left( U_{n-1}^{(1)}, \ldots,  U_{n-1}^{(n-1)}\right),
\ee
where $( U_{n-1}^{(i)})_{1\le i\le n-1}$ 
are the order statistics of a sample 
$(s+(1-s)U_i)_{1\le i\le n-1}$, with 
$(U_i)_{1\le i\le n-1}$ i.i.d.  uniform $[0,1]$. Thus for $0<s< 1\le x$
\ben
\PP(s+(1-s)U_1\le x^{-\alpha}) =
\PP\bigg(U_1\le \frac{x^{-\alpha}-s}{1-s}\bigg)
=\frac{x^{-\alpha}-s}{1-s}.
\een
This  equals $\PP(J_1(s^{1/\alpha})\le x^{-\alpha})$ as calculated from \eqref{Jdis0} so we get the required representation in \eqref{top-}.
When $r=0$, \eqref{df} remains true with  $( U_{n-1}^{(i)})_{1\le i\le n-1}$ the order statistics of $(U_i)_{1\le i\le n-1}$ i.i.d.  uniform $[0,1]$, and since $\PP(U_1\le x^{-\alpha})=x^{-\alpha}=\PP(L_1>x)$, with $L_1$ as in \eqref{Ldis}, we get  \eqref{top-L}.

Next  suppose $\pibar^+$ is continuous at $xq_t$, $x>0$, and $xq_t$ 
  is a point of decrease of $\pibar^+$.
  Hold $t>0$ fixed.
    The continuous case analogue of the ratio in  \eqref{rat} is
\be\label{x11}
\lim_{\veps\dto 0}
\frac{\PP\big(\Delta X_t^{(r+k)}>x_k\Delta X_t^{(r+n)},\, 1\le k\le n-1,\, xq_t-\veps <\Delta X_t^{(r+n)}\le xq_t+\veps\big)}
{\PP\big(xq_t-\veps<\Delta X_t^{(r+n)}\le xq_t+\veps\big)}.
\ee
Letting $a_t(x,\veps):= t\pibar^+(xq_t+\veps)
<b_t(x,\veps):= t\pibar^+(xq_t-\veps)$ for 
 $\veps\in(0,xq_t)$, the denominator in \eqref{x11} is
 \be\label{x2}
 \int_{a_t(x,\veps)}^{b_t(x,\veps)} f_{r+n}(y)\rmd y
 = (b_t(x,\veps)-a_t(x,\veps)) f_{r+n}(\xi_t(x,\veps)),
 \ee
for some  $\xi_t(x,\veps)\in [a_t(x,\veps),b_t(x,\veps))$.
Note that the righthand side of \eqref{x2} is positive since 
 $xq_t$   is a point of decrease of $\pibar^+$.
 Let $c_t(x_k,x,\veps):= t\pibar^+(x_k(xq_t-\veps))$. 
 By a similar calculation as in \eqref{x1}
 (but noting the inequalities in \eqref{x11} as opposed to the equality in \eqref{rat}),
 the numerator in \eqref{x11} is not greater than 
 \bean %%\label{x12}
 &&
\int_{a_t(x,\veps)}^{b_t(x,\veps)} 
\PP\bigg(\frac{\Gamma_{r+k}}{\Gamma_{r+n}}< \frac{c_t(x_k,x,\veps)}{y},\, 1\le k\le n-1\bigg)     f_{r+n}(y)\rmd y\cr
&&\cr
&&=
(b_t(x,\veps)-a_t(x,\veps)) f_{r+n}(\eta_t(x,\veps))
\PP\bigg(\frac{\Gamma_{r+k}}{\Gamma_{r+n}}< \frac{c_t(x_k,x,\veps)}{\eta_t(x,\veps)},\, 1\le k\le n-1\bigg) 
   \eean
where $\eta_t(x,\veps)\in [a_t(x,\veps),b_t(x,\veps))$.
Letting $\veps\dto 0$ we find an upper bound of the form
 \bean %%\label{x3}
 &&
\PP\bigg(\Delta X_t^{(r+k)}>x_k\Delta X_t^{(r+n)},\, 1\le k\le n-1\big| \Delta X_t^{(r+n)}=xq_t\bigg)\cr
&&\le 
\PP\bigg(\frac{\Gamma_{r+k}}{\Gamma_{r+n}}< \frac{t\pibar^+(x_kxq_t-)}{t\pibar^+(xq_t)},\, 1\le k\le n-1\bigg) 
\eean
at points $t>0$, $x>0$, such that $xq_t$ is a point of decrease of $\pibar^+$.
Similarly we get a lower bound 
%%for the conditional probability on the lefthand side 
with  $t\pibar^+(x_kxq_t)$ replacing $t\pibar^+(x_kxq_t-)$.
Then as $t\dto 0$, on account of the regular variation of $\pibar^+$, both  bounds approach the expression in \eqref{de}, which can be re-expressed in terms of the $J_i$ and $L_i$, as shown. 
Having reached the same limit in both cases, we have proved  Proposition \ref{large_j}.
\halmos

\section{Proofs for Section \ref{s5}}\label{ps5}
%%: Convergence of Trimmed L\'evys}
Throughout this section $X$ will be a  L\'evy process in the domain of attraction at 0 of a non-normal stable  random variable. Thus the  L\'evy tail $\pibar$ is  regularly varying of index $-\alpha$, $\alpha\in(0,2)$, at 0, and the balance condition \eqref{bal} holds at 0. Since $a_+>0$ in \eqref{bal},  also
 $\pibar^+\in RV_0(-\alpha)$ at 0.
% (or at $\infty$).  We identify $\Delta X_t=X_t-X_{t-}$ with the Poisson points   $\Delta_t$ in  Section \ref{sec: PPP}.
%
% Since there is no negative component to $\Pi$ we simply write 
% $\pibar(x)$ and $\pibarinv(x)$  for $\pibar^+(x)$ and 
% $\pibarpinv(x)$ from now on.
% The trimmed stable  random variable  ${}^{(r)}S_1$ and the stable process jumps are as in Section \ref{??}.

\medskip\noindent {\bf Proof of  Theorem \ref{stable}:}\
(i)\ 
Take $r\in\N_0$, $n\in\N$, and choose $x_1\ge \cdots \ge x_{n-1}\ge 1$, $x_n=1$, $\theta_k\in\R$, $1\le k\le n$,
and $v>0$.
For shorthand, write $M_t^{(r+n)}$ for $\rho_X(\Delta X_t^{(r+n)})$.
We proceed by finding the limit as $t\dto 0$ of the conditional  characteristic function 
\begin{align}\label{6a}
&
\EE\bigg(\exp\bigg(\rmi \sum_{k=1}^n 
\frac{\theta_k({}^{(r)}X_t-tM_t^{(r+n)})}{\Delta X_t^{(r+k)}}\bigg)
\bigg|\frac{\Delta X_t^{(r+k)}}{\Delta X_t^{(r+n)}}
=x_k, 1\le k<n,\frac{\Delta X_t^{(r+n)}}{\pibarpinv(1/t)}=v^{-1/\alpha}
\bigg)\cr
&=
\EE\bigg(\exp\bigg(\rmi\sum_{k=1}^n \frac{\theta_k}{x_{k}}
\frac{({}^{(r)}X_t-tM_t^{(r+n)})}{\Delta X_t^{(r+n)}}\bigg)\cr
&
\hskip1.5cm 
\bigg|\frac{\Delta X_t^{(r+k)}}{\Delta X_t^{(r+n)}}
=x_k, 1\le k\le n-1,\, \frac{\Delta X_t^{(r+n)}}{\pibarpinv(1/t)}
=v^{-1/\alpha}\bigg).
\end{align}
Decompose ${}^{(r)}X_t$ as follows:
\be\label{3c-}
 \frac{{}^{(r)}X_t}{\Delta X_t^{(r+n)}} 
= \sum_{k=1}^{n}
\frac{\Delta X_t^{(r+k)}}{\Delta X_t^{(r+n)}} + 
\frac{{}^{(r+n)}X_t}{\Delta X_t^{(r+n)}},
 \ee
and recall the definitions of $x_{n+}$ and $\wt\theta_{n+}$ in \eqref{thep}.
%%where  by convention we set $\sum_{k=1}^0=0$.
%For brevity define 
%\be\label{thep}
%x_{n+}=\sum_{k=1}^nx_k \quad {\rm and}\quad
%\wt\theta_{n+}=\wt\theta_{n+}(x_1,\ldots,x_n):=  \sum_{k=1}^n \frac{\theta_k}{x_k}.
%\ee
Given the conditioning in \eqref{6a}, 
the first component on the RHS of \eqref{3c-} equals $\sum_{k=1}^nx_{k}=x_{n+}$, so we can 
%%For the second component, recall $x_n=1$ and apply Corollary \ref{pr}  
%with $r$ replaced by $r+n$ to replace the conditioning on $\Delta X_t^{(r+k)}$, $1\le k\le n$, in \eqref{6a},
% by conditioning on 
%$\Delta X_t^{(r+n)}$.
%Then to 
write the RHS of \eqref{6a} as % $e^{\rmi\wt\theta_{n+} x_{n+}}$ times
\be\label{7a}
e^{\rmi\wt\theta_{n+} x_{n+}}\times
\EE\bigg(\exp\bigg(\rmi\wt\theta_{n+}
\frac{{}^{(r+n)}X_t-tM_t^{(r+n)}}{\Delta X_t^{(r+n)}}\bigg)
\bigg|\frac{\Delta X_t^{(r+k)}}{\pibarpinv(1/t)}=x_kv^{-1/\alpha},\, 1\le k\le n\bigg)
\ee
(recall $x_{n}=1$). 
Then by \eqref{B9} with $\theta$ replaced by $\wt\theta_{n+}$, 
%$u_r$ replaced by $\pibarinv(v/t)$, 
$y_k$ replaced by $y_k(t):= x_kv^{-1/\alpha}\pibarpinv(1/t)$,  and $\sigma^2=0$, 
%and $\rmd_X=0$, 
the expression in \eqref{7a} equals
\bea\label{8}
&&
e^{\rmi\wt\theta_{n+} x_{n+}}\times
\exp\bigg(\int_{(-\infty,1)}\big(e^{\rmi\wt\theta_{n+}x}-1
-\rmi\wt\theta_{n+}x{\bf 1}_{\{|x|\le 1\}}\big)
t\Pi\big(v^{-1/\alpha}\pibarpinv(1/t)\rmd x\big)\bigg)\cr
&&\hskip3.8cm 
\times K_{r+1,r+n}(\wt\theta_{n+}/y_{r+n}(t),\, t,\, y_{r+1}(t),\ldots, y_{r+n}(t))
\eea
%%Now replace $x_{r+n}$ by $v^{-1/\alpha}\pibarpinv(1/t)$, where $v>0$.
(again, recall $x_{n}=1$). 
By \eqref{1b}, we have
$t\pibar^+(v^{-1/\alpha}\pibarpinv(1/t))\to v$,
and hence 
$t\pibar^+(v^{-1/\alpha}\pibarpinv(1/t)\rmd x)\to v\Lambda(\rmd x)$, $x>0$, 
vaguely, as $t\dto 0$.
The limit of the second factor in \eqref{8} can then be found straightforwardly using integration by parts and applying \eqref{bal} and \eqref{1b}.

The term containing $K$ in \eqref{8} is negligible here, as follows.
Note that  $\pibar^+$ in $RV_0(-\alpha)$ implies
$\Delta\pibar^+(x):=\pibar^+(x-)-\pibar^+(x)=o(\pibar^+(x))$ as $x\dto 0$.
Recall  $K_{r,r+n}$ is defined in  \eqref{Krdef}, and $\kappa(y,v)=\pibar^+(y-)-v$. 
%%When $x_{r+n}/t\in(a_r,b_r)= (\pibar^+(x_{r+n}),\pibar^+(x_{r+n}-))$
Substituting $y_{r+n}(t)=v^{-1/\alpha}\pibarpinv(1/t)$ for $y$ gives
\bea\label{kau}
t\kappa(y_{r+n}(t),v/t)&=& 
t\pibar^+(v^{-1/\alpha} \pibarpinv(1/t)-)-v\cr
&=&
t\pibar^+(v^{-1/\alpha} \pibarpinv(1/t))-v+t\Delta\big(\pibar^+(v^{-1/\alpha} \pibarpinv(1/t))\big)\cr
&=&
t\pibar^+(v^{-1/\alpha} \pibarpinv(1/t))-v+o\big(t\pibar^+(v^{-1/\alpha} \pibarpinv(1/t))\big)\cr
&\to&
v-v=0,\ {\rm as}\ t\dto 0.
\eea
Furthermore, $G_t^{(y_{r+n}(t),\Gamma_{r+n}/t)}/y_{r+n}(t)$ has the distribution of $Y_{t\kappa(y_{r+n}(t),\Gamma_{r+n}/t)}$ and hence tends to 0 in probability  when $t\dto 0$.
So we can ignore the $K$ term in \eqref{8}.

We conclude that the expression in  \eqref{8}  tends  as $t\dto 0$ to 
\be\label{10}
e^{\rmi\wt\theta_{n+} x_{n+}}\times
\exp\bigg(v\int_{(-\infty,1)}\big(e^{\rmi\wt\theta_{n+}x}-1-\rmi\wt\theta_{n+}x{\bf 1}_{\{|x|\le 1\}}\big)
\Lambda(\rmd x)\bigg).
\ee
%%where $x_{n+}:= \sum_{k=1}^nx_k$. 
Thus, by \eqref{3c-}, to find  the limit as $t\dto 0$ of
\ben  %%\label{6a}
\EE\exp\bigg(\rmi \sum_{k=1}^n 
\frac{\theta_k( {}^{(r)}X_t-tM_t^{(r+n)})}{\Delta X_t^{(r+k)}}\bigg),
\een
we multiply \eqref{10} by the limit, 
as $t\dto 0$ through points $v$ such that $v^{-1/\alpha}\pibarpinv(1/t)$ is a point of decrease of $\pibar^+$, of
\bean  %\label{pra}
&&
\PP\bigg(
\frac{\Delta X_t^{(r+k)}}{\Delta X_t^{(r+n)}}\in\rmd x_k,\, 1\le k\le n-1
\bigg|\,
\frac{\Delta X_t^{(r+n)}}{ \pibarpinv(1/t)}=v^{-1/\alpha}\bigg)\cr
&&
\hskip7cm   \times 
\PP\bigg(\frac{\Delta X_t^{(r+n)}}{ \pibarpinv(1/t)}\in \rmd(v^{-1/\alpha})\bigg),
\eean
and then integrate over $v$ and the $x_k$.\footnote{We use the result:
$\int f_t(\omega)P_t(\rmd\omega)\to \int f(\omega)\PP(\rmd\omega)$ 
when $P_t\toweak P$ are probability measures and $f_t\to f$, $f$ continuous, $|f|\le 1$.
In \eqref{6a}, the $f_t$ are characteristic functions and the limit distribution $P$ in \eqref{25} is continuous in all its variables.}

From  \eqref{top-} %%with $x_0=1$ %%in Theorem \ref{pre} 
when $r\in\N$ and from \eqref{top-L} when $r=0$ 
we see that the limit of the conditional probability depends only on 
the $J_{n-1}^{(k)}$ or $L_{n-1}^{(k)}$ and $B_{r,n}$, 
and not on $v$, 
while by \eqref{1b}
\bean
\PP\big(\Delta X_t^{(r+n)}> v^{-1/\alpha} \pibarpinv(1/t)\big)
&=&
\PP\big(\pibarpinv(\Gamma_{r+n}/t) >v^{-1/\alpha} \pibarpinv(1/t)\big)\cr
&=&
\PP\big(\Gamma_{r+n}< t\pibar^+(v^{-1/\alpha} \pibarpinv(1/t))\big)\cr
&\to&
\PP\big(\Gamma_{r+n}\le v\big), \ {\rm as}\ t\dto 0.
\eean
%%%%NB:  \pibarinv(v) > x iff v < \pibar(x)
Putting the RHS of \eqref{top-} or  \eqref{top-L} together
with the expression in \eqref{10}
we can write the limiting characteristic function of 
%\ben
%\left( \frac{{}^{(r)}X_t}{\Delta X_t^{(r+1)}}, \ldots, 
%\frac{{}^{(r)}X_t} {\Delta X_t^{(r+n)}}\right) 
%\een
the $n$-vector on the LHS of \eqref{liPD} 
as %%(recall $x_n=1$)
\begin{align} \label{25}
&
\int_{\bfx^\uparrow\ge 1}
e^{\rmi\wt\theta_{n+}x_{n+}}
\int_0^\infty \exp\bigg(v\int_{-\infty}^1
\big(e^{\rmi\wt\theta_{n+}x}-1-\rmi\wt\theta_{n+}x{\bf 1}_{\{|x|\le 1\}}\big)
\Lambda(\rmd x)\bigg)\PP\left(\Gamma_{r+n}\in\rmd v\right)\nonumber\\
&\hskip6cm 
\times 
\PP\big(J_{n-1}^{(k)}\big({\rm B}_{r,n}^{1/\alpha}\big)\in\rmd x_k, 1\le k\le n-1\big)
\end{align}
when $r\in\N$, and with each
$J_{n-1}^{(k)}({\rm B}_{r,n}^{1/\alpha})$
 replaced by $L_{n-1}^{(k)}$ when $r=0$.
%  (We replaced the integration variables $x_{r+1}, \ldots,x_{r+n-1}$ by
%  $x_1,\ldots, x_{n-1}$, and 
Recall that $\int_{\bfx^\uparrow\ge 1}$ denotes  
integration over the region $\{x_1\ge x_2\ge \cdots \ge x_{n-1}\ge 1\}$.

Note that, with $\Lambda$ defined as in \eqref{LM}, $\lambar(x)\in RV_0(-\alpha)$ and $\lambar^+$ and $\lambar^-$ satisfy \eqref{bal}.
So exactly the same calculation\footnote{This easy correspondence is the reason for adopting the nonstandard centering in \eqref{ces}.}
 with
$\Delta S_1^{(k)}$ replacing $\Delta X_t^{(k)}$, for $r+1\le k\le r+n$, 
and $\Lambda$ replacing $\Pi$ (and no limit on $t$ is necessary), shows that the  characteristic function of the vector of  stable ratios on the RHS of \eqref{liPD} equals 
\eqref{25} when $r\in\N$ or the corresponding version 
when $r=0$.

(ii)\
To derive \eqref{38} and  the corresponding version 
when $r=0$, observe that
%$\exp\left(v\int_{(-\infty,1)}\left(e^{\rmi\theta x}-1-\rmi\theta x{\bf 1}_{\{|x|\le 1\}}\right) \Lambda(\rmd x)\right)$, $\theta\in\R$, 
the exponent inside the integral in  \eqref{25} 
is the characteristic function of a L\'evy process $(W_v)_{v\ge 0}$ having L\'evy triplet $(0,0,\Lambda(\rmd x){\bf 1}_{(-\infty,1)})$, that is, of
a Stable$(\alpha)$ process with jumps truncated below 1.
So the integral with respect to $v$ in \eqref{25} is
\ben %%\label{21}
%%\int_{\bfx^\uparrow\ge 1}e^{\rmi\wt\theta_{n+}x_{n+}}
\int_{v>0} \EE\big(e^{\rmi\wt\theta_{n+}  W_v}\big)
\PP\left(\Gamma_{r+n}\in\rmd v\right)
%%=\int_{\bfx^\uparrow\ge 1}e^{\rmi\wt\theta_{n+}x_{n+}}
=
\EE\big(e^{\rmi\wt\theta_{n+} W_{\Gamma_{r+n}}}\big),
\een
%In the second factor in \eqref{25},  the integration over $x_0$
%can be written as an expectation with respect to an  independent ${\rm B}_{r,n}$  random variable. 
and thus we obtain \eqref{38} when $r\in\N$
and the corresponding version when $r=0$ with the $J_i$ replaced by $L_i$.

When $r\in\N$ and $n=2,3,\ldots$, the alternative representation in \eqref{39} is obtained by evaluating the $\rmd v$ integral in \eqref{25}, resulting in
(recall $\psi(\cdot)$ defined in \eqref{3a01}):
\begin{align}\label{27}
&
\EE\exp\bigg(\rmi \sum_{k=1}^n 
\frac{\theta_k\big({}^{(r)}X_t-t\rho_X(\Delta X_t^{(r+n)})\big)}{\Delta X_t^{(r+n)}}\bigg)
\nonumber \\
&
\to
\EE\exp\bigg(\rmi \sum_{k=1}^n 
\frac{\theta_k\big({}^{(r)}S_1-\rho_S(\Delta S_1^{(r+n)})\big)}{\Delta S_1^{(r+n)}}\bigg)
\nonumber \\
&
=\int_{\bfx^\uparrow\ge 1}
e^{\rmi\wt\theta_{n+}x_{n+}}
\int_{v>0}
\frac{v^{r+n-1}e^{-v\left(1-\psi(\wt\theta_{n+})\right)}}{\Gamma(r+n)}\rmd v\
\nonumber \\
& \hspace{1in}
\times \PP\big(J_{n-1}^{(k)}({\rm B}^{1/\alpha}_{r,n})\in\rmd x_k, 1\le k\le n-1\big),
\end{align}
 equal to the expression in \eqref{39}.
%  (Note that $\psi(0)=0$, so $|\psi(\theta)|<1$ for $|\theta|\le$ some $ \theta_0>0$.)
% &&\cr
%&&=
%\int_{\bfx^\uparrow\ge 1}
%%%\int_{(x_1,\ldots, x_{n-1})\in\R^{n-1}}
%\frac{e^{\rmi\wt\theta_{n+}x_{n+}}}
%{\left(1-\psi(\wt\theta_{n+})\right)^{r+n}}
%  \PP\left(J_{n-1}^{(k)}({\rm B}^{1/\alpha}_{r,n})\in\rmd x_k, 1\le k\le n-1\right).
%\eea
%\bea\label{22}
%&&
%\lim_{t\dto 0}
%\EE\left(\exp\left(\rmi \sum_{k=1}^n 
%\frac{\theta_k {}^{(r)}X_t}{\Delta X_t^{(r+k)}}\right)\right)=
%\EE\left(\exp\left(\rmi \sum_{k=1}^n 
%\frac{\theta_k {}^{(r)}S_1}{\Delta S_1^{(r+k)}}\right)\right)
%\cr
%&&\cr
%&&=
%\int_{\bfx^\uparrow\ge 1}
%\EE\left(e^{\rmi\wt\theta_{n+}(x_{n+}+W_{\Gamma_{r+n}})}\right) 
%\PP\left(J_{n-1}^{(k)}({\rm B}^{1/\alpha}_{r,n})\in\rmd x_k, 1\le k\le n-1\right),\cr
%&&
%\eea
%where $\Gamma_{r+n}$ is  independent of $W$, ${\rm B}_{r,n}$ is  independent of the $J_{n-1}^{(k)}$, and $ \EE_{B_{r,n}}$ denotes expectation with respect to $B_{r,n}$.
 When $r\in\N_0$, $n=1$, similar working shows that 
\eqref{25} can be replaced by
\bea\label{26}
&&
\lim_{t\dto 0}
\EE\big(e^{\rmi \theta \big({}^{(r)}X_t-t\rho_X(\Delta X_t^{(r+1)})\big)/
\Delta X_t^{(r+1)}}\big)\cr
&&
=e^{\rmi\theta}\times
\int_{v>0}\exp\bigg(v\int_{(-\infty,1)}\big(e^{\rmi\theta x}-1-\rmi\theta x{\bf 1}_{\{|x|\le 1\}}\big)
\Lambda(\rmd x)\bigg)\PP\left(\Gamma_{r+1}\in\rmd v\right)\cr
&&\cr
&&=
e^{\rmi\theta}\times\EE\big(e^{\rmi\theta W_{\Gamma_{r+1}}}\big),\,
\theta\in\R.
\eea
%%gives the limiting characteristic function of 
%%${}^{(r)}X_t/\Delta X_t^{(r+1)}$, for $t\dto 0$.  
 
%%   \medskip\noindent {\bf Proof of  Corollary \ref{prx}:}\

(iii)\
Finally, to prove   \eqref{36} when $r\in\N$, set $\theta_1=\cdots=\theta_{n-1}=0$, $\theta_n=\theta$ (so $\wt\theta_{n+}=\theta$ and, recall, $x_{n+}=x_1+\cdots +x_{n-1}+1$)
in \eqref{27} to get the characteristic function of the RHS of \eqref{liPD0} equal to 
\bean %%\label{36a}
&&
\int_{\bfx^\uparrow\ge 1}
\frac{e^{\rmi\theta x_{n+}}}{\left(1-\psi(\theta)\right)^{r+n}}
 \PP\big(J_{n-1}^{(k)}({\rm B}^{1/\alpha}_{r,n})\in\rmd x_k, 1\le k\le n-1\big)\cr
 &&\cr
 &&=
 \frac{e^{\rmi\theta }}{\left(1-\psi(\theta)\right)^{r+n}}
 \int_{0<u<1}\EE
 \exp\bigg(\rmi\theta\sum_{k=1}^{n-1}J_{n-1}^{(k)}(u)\bigg)
 \PP\big({\rm B}^{1/\alpha}_{r,n}\in\rmd u\big)\cr
 &&\cr
 &&=
  \frac{e^{\rmi\theta }}{\left(1-\psi(\theta)\right)^{r+n}}
 \int_{0<u<1} \big(\EE e^{\rmi\theta J_1(u)}\big)^{n-1}
 \PP\big({\rm B}^{1/\alpha}_{r,n}\in\rmd u\big)\cr
 &&\cr
 &&=
   \frac{e^{\rmi\theta }}{\left(1-\psi(\theta)\right)^{r+n}}
   \EE\big( \phi^{n-1}(\theta, {\rm B}^{1/\alpha}_{r,n})\big),
\eean
where $\phi(\theta,u)=\EE e^{\rmi\theta J_1(u)}$ as in \eqref{phiJ},
with $|\psi(\theta)|<1$ when  $|\theta|\le\theta_0$.
Similarly, \eqref{26} can alternatively be written  as $e^{\rmi\theta}$ times  the expression in \eqref{31}.
The $r=0$ case follows as before. \halmos
%\be\label{28}
%\lim_{t\dto 0}
%\EE\left(e^{\rmi \theta {}^{(r)}X_t/
%\Delta X_t^{(r+1)}}\right)=
%\EE\left(e^{\rmi \theta {}^{(r)}S_1/
%\Delta S_1^{(r+1)}}\right)=
%e^{\rmi\theta}\EE\left(e^{\rmi\theta W_{\Gamma_{r}}}\right)
%=
%\frac{e^{\rmi\theta}}
%{\left(1-\psi(\theta)\right)^{r+1}}.
%\ee
%%Hence the corollary. \halmos

\medskip\noindent {\bf Proof of Theorem \ref{prin}:}\
In this proof $X$ is a driftless subordinator whose L\'evy tail measure is in $RV_0(-\alpha)$, $0<\alpha<1$.
From \eqref{2333} we obtain the Laplace transform
\begin{align}\label{B13}
&\EE\exp\bigg(-\lambda\frac{{}^{(r+n)}X_t}{\Delta X_t^{(r)}}\bigg)=
\int_{y>0}\int_{w>y} e^{-t\int_{(0,a)}\left(1-e^{-\lambda x/b}\right)\Pi(\rmd x)
-t\kappa(w/t)(1-e^{-\lambda a/b})}
\nonumber \\
&
\hskip8cm 
\times \PP\left(\Gamma_{r}\in \rmd y, \Gamma_{r+n}\in \rmd w\right),
\end{align}
where $\lambda>0$ and for brevity
\ben
a=a(w,t):= \pibarinv(w/t)\le b=b(y,t):=\pibarinv(y/t),\ t>0,\  w>y>0
\een
(we can write $\pibar$ and $\pibarinv$ for $\pibar^+$ and $\pibarpinv$ in  \eqref{B11}).
We derive an upper bound for the exponent in \eqref{B13} as follows.
%Now $\Delta \pibar(x)=0\left(\pibar(x)\right)$ as $x\dto 0$ (cf. \eqref{kau}) so there is a $\delta>0$ such that 
% $\Delta \pibar(x)\le \pibar(x)$ for $0<x\le \delta$.
% We may choose $z_1$ larger if necessary so that $\pibarinv(z_1)\le \delta$.
% Then, as in \eqref{kau}, $y/t\ge z_1$ implies $\pibarpinv(y/t)\le \pibarinv(z_1)\le \delta$ and (see \eqref{B13}) 
% \ben
% t\kappa(y/t)\big(1-e^{-\lambda}\big) \le t\lambda\Delta \pibar(\pibarinv(y/t))\le \lambda y.
% \een
Keep $0<t\le t_0$ for a fixed $t_0>0$, throughout.

First, the integral in the exponent of \eqref{B13} is
\bea \label{B21}
t\int_{(0,a)}\big(1-e^{-\lambda x/b}\big)\Pi(\rmd x)
&\le &
t(\lambda/b)\int_0^a e^{-\lambda x/b}\pibar(x)\rmd x\quad {\rm (integrate\ by\ parts)}\nonumber \\
&=&
 t\lambda\int_0^{a/b}e^{-\lambda x}\pibar(bx)\rmd x.
\eea
Now $a(w,t)\to \pibarinv(+\infty)=0$ as $w\to\infty$ or $t\dto 0$, and 
 $b(y,t)\to \pibarinv(+\infty)=0$ as $y\to\infty$ or $t\dto 0$.
To compare the magnitudes of $a$ and $b$ we use the Potter bounds 
%Bingham, N.H., Goldie, C.M. and Teugels, J.L. (1987)
%Regular Variation. Encyclopedia of Mathematics and its Applications,27, Cambridge University Press, Cambridge. 
 (Bingham, Goldie and Teugels \citeyearpar[p.25]{BGT87}). 
Since $\pibar\in RV_0(-\alpha)$ with $0<\alpha<1$, given $\eta>0$ 
%with  $\eta<1/\alpha$
there are constants $c>0$ and $z_0=z_0(\eta)>0$ such that
\be\label{B19}
\frac{\pibar(\mu z)}{\pibar(z)} \le c\mu^{-\alpha-\eta}\ {\rm for\ all}\ \mu\in(0,1),\ z\in(0,z_0);
\ee
and since $\pibarinv\in RV_\infty(-1/\alpha)$ we also have
\be\label{B20}
\frac{\pibarinv(\mu z)}{\pibarinv(z)} \le c\mu^{-1/\alpha+\eta}\ {\rm for\ all}\ \mu>1,\ z>1/z_0
\ee
(where $c$ and $z_0$ may be chosen the same in both cases, and  $\eta<1/\alpha$).
Thus for $0<x\le a/b\le 1$ and $0<b\le z_0$, using \eqref{B19},
\ben
t\pibar(bx)\le ct x^{-\alpha-\eta}\pibar(b)
= ct x^{-\alpha-\eta}\pibar(\pibarinv(y/t))\le 
cyx^{-\alpha-\eta},
\een
and we have $b\le z_0$ if $\pibarinv(y/t)\le z_0$, i.e., if $y/t\ge \pibar(z_0)$. For $w>y$ and $y/t\ge 1/z_0$, using \eqref{B20}, 
\be\label{a/b}
\frac{a}{b}=\frac{\pibarinv(w/t)}{\pibarinv(y/t)}
=\frac{\pibarinv((w/y)(y/t))}{\pibarinv(y/t)}
\le c\left(\frac{w}{y}\right)^{-1/\alpha+\eta}
=  c\left(\frac{y}{w}\right)^{1/\alpha-\eta}.
\ee
Now keep $y/t\ge z_1:= \pibar(z_0)\vee(1/z_0)$ and
$0<\eta<1-\alpha$ (so also $\eta<1/\alpha$). Then by \eqref{B21}
\bea\label{I1}
&& 
t\int_{(0,a)}\big(1-e^{-\lambda x/b}\big)\Pi(\rmd x)
\le t\lambda\int_0^{a/b}e^{-\lambda x}\pibar(bx)\rmd x\cr
&&\le 
  c\lambda y\int_0^{a/b}x^{-\alpha-\eta}\rmd x
=  \frac{c\lambda y}{1-\alpha-\eta} \left(\frac{a}{b}\right)^{1-\alpha-\eta}\cr
&&\le  c'\lambda y\left(\frac{y}{w}\right)^\beta=: 
\lambda  g_1(w,y),
 \eea
where 
 $c':= c^{2-\alpha-\eta}/(1-\alpha-\eta)>0$ and 
$\beta:= (1-\alpha-\eta)(1/\alpha-\eta)>0$.

Alternatively, when  $y/t<z_1$, we have $b=\pibarinv(y/t)\ge \pibarinv(z_1)$, while $t\le t_0$ implies  $a=\pibarinv(w/t)\le \pibarinv(w/t_0)$.  Then
\bea\label{B23}
t\int_{(0,a)}\big(1-e^{-\lambda x/b}\big)\Pi(\rmd x)
&\le& 
t(\lambda/b)\int_{(0,a)}x\Pi(\rmd x)\cr
&\le&
t_0(\lambda/\pibarinv(z_1))\int_{(0,\pibarinv(w/t_0))}x\Pi(\rmd x)\cr
&=:&\lambda  g_2(w).
 \eea
 
For the term containing $\kappa$ in \eqref{B13},  we have, for all $x\in(0,z_0)$, 
\ben
\Delta \pibar(x)=\pibar(x-)-\pibar(x) \le 
\pibar(x/2)-\pibar(x)\le 2^{\alpha+\eta} c\pibar(x)
\een
by \eqref{B19}.
Thus for all $t>0$ and $w>y>0$, using \eqref{kau}, 
\bean
t\kappa(w/t)\le t\Delta  \pibar(\pibarinv(w/t)) 
\le
2^{\alpha+\eta} c t\pibar(\pibarinv(w/t)) 
\le 2^{\alpha+\eta} cw,
\eean 
because we kept $y/t>\pibar(z_0)$, 
as a consequence of which 
$\pibarinv(w/t)\le \pibarinv(y/t)\le z_0$. 
%%and so $t\kappa(w/t)\le cw$  for all $t>0$ and $w>0$.
Then $t\kappa(w/t)(1-e^{-\lambda  a/b})
\le c 2^{\alpha+\eta}  w\lambda a/b$.
When $w>y$ and $y/t\ge 1/z_0$, $t\kappa(w/t)(1-e^{-\lambda  a/b})
\le 2^{\alpha+\eta}  c^2\lambda w(y/w)^{1/\alpha-\eta}$ by \eqref{a/b}.
 When $y/t<z_1$, so  $b\ge \pibarinv(z_1)$, 
 $t\kappa(w/t)(1-e^{-\lambda  a/b})
\le  2^{\alpha+\eta} cw\lambda \pibarinv(w/t_0)/\pibarinv(z_1)$.
  So an overall upper bound  for the term containing $\kappa$ 
in \eqref{B13} is
\begin{align}\label{B22}
& t\kappa(w/t)(1-e^{-\lambda  a/b}) \nonumber \\
&\le 
\lambda  g_3(w,y):=2^{\alpha+\eta} 
\max\big(c^2\lambda w(y/w)^{1/\alpha-\eta}, \,cw\lambda \pibarinv(w/t_0)/\pibarinv(z_1)\big).
%\le
%cw\lambda \pibarinv(w/t_0)=:
\end{align}
    
  Combine \eqref{I1}--\eqref{B22} to get an upper bound for the negative of the  exponent in \eqref{B13} of the form 
\ben
\lambda g(w,y):=\lambda\big( \max(g_1(w,y),g_2(w))+ g_3(w,y)\big).
\een
So, for all $0<t\le t_0$, $n\in\N$, 
  \bea\label{B16}
\EE\exp\left(-\lambda\frac{{}^{(r+n)}X_t}{\Delta X_t^{(r)}}\right)
&\ge&
\int_{y>0}\int_{w>y}e^{-t_0\lambda   g(w,y)}
\PP\left(\Gamma_{r}\in \rmd y, \Gamma_{r+n}\in \rmd w\right)\cr
&=&
\EE\big(e^{-t_0\lambda   g(\Gamma_{r+n},\Gamma_{r})}\big).
\eea
Now when $w\to\infty$, $g_1(w,y)\to0$ 
 for each $y>0$ (see \eqref{I1}) and
  $g_2(w)\to 0$ as $w\to\infty$ 
   because  $\pibarinv(w)\to 0$ as $w\to\infty$ (see \eqref{B23});
  while    $g_3(w,y)\to0$  for each $y>0$ 
because $\pibarinv\in RV_\infty(-1/\alpha)$ and $0<\alpha<1$ (see \eqref{B22}).

Finally, since $\Gamma_{r+n}\topr\infty$ as $n\to\infty$ for each $r\in\N$, we can let $n\to\infty$ and use Fatou's lemma in \eqref{B16} to see that 
\bean\label{B25}
\EE\exp\left(-\lambda\frac{{}^{(r+n)}X_t}{\Delta X_t^{(r)}}\right)
&\to& 1,\ {\rm as}\ n\to\infty,
\eean
for each $r\in\N$, uniformly in $\lambda>0$ and $t\in(0,t_0]$.
We deduce convergence  in probability in \eqref{B15} uniformly in $t\in(0,t_0]$ from this.
%%, then since the LHS of \eqref{B15} is monotone in $n$, we get the a.s. convergence.\halmos

%\medskip\noindent {\bf Acknowledgements.}\
%We are grateful for helpful feedback from 
%P\'eter  Kevei and David Mason, and from a referee.

% 
% 
 \renewcommand{\bibfont}{\small}
 \bibliography{Library_Levy_Feb2018}
 \bibliographystyle{newapa}

\end{document}